\magnification\magstep1
\baselineskip15pt

\newread\AUX\immediate\openin\AUX=\jobname.aux
\newcount\relFnno
\def\ref#1{\expandafter\edef\csname#1\endcsname}
\ifeof\AUX\immediate\write16{\jobname.aux gibt es nicht!}\else
\input \jobname.aux
\fi\immediate\closein\AUX

\newcount\itemcount
\def\resetitem{\global\itemcount0}\resetitem
\newcount\Itemcount
\Itemcount0
\newcount\maxItemcount
\maxItemcount=0

\def\FILTER\fam\itfam\tenit#1){#1}

\def\Item#1{\global\advance\itemcount1
\edef\TEXT{{\it\romannumeral\itemcount)}}%
\ifx?#1?\relax\else
\ifnum#1>\maxItemcount\global\maxItemcount=#1\fi
\expandafter\ifx\csname I#1\endcsname\relax\else
\edef\testA{\csname I#1\endcsname}
\expandafter\expandafter\def\expandafter\testB\testA
\edef\testC{\expandafter\FILTER\testB}
\edef\testD{\csname0\testC0\endcsname}\fi
\edef\testE{\csname0\romannumeral\itemcount0\endcsname}
\ifx\testD\testE\relax\else
\immediate\write16{I#1 hat sich geaendert!}\fi
\expandwrite\AUX{\neverexpand\ref{I#1}{\TEXT}}\fi
\item{\TEXT}}

\def\today{\number\day.~\ifcase\month\or
  Januar\or Februar\or M{\"a}rz\or April\or Mai\or Juni\or
  Juli\or August\or September\or Oktober\or November\or Dezember\fi
  \space\number\year}
\font\sevenex=cmex7
\scriptfont3=\sevenex
\font\fiveex=cmex10 scaled 500
\scriptscriptfont3=\fiveex

\def\phi{\varphi}
\def\epsilon{\varepsilon}
\def\theta{\vartheta}
\def\uauf{\lower1.7pt\hbox to 3pt{%
\vbox{\offinterlineskip
\hbox{\vbox to 8.5pt{\leaders\vrule width0.2pt\vfill}%
\kern-.3pt\hbox{\lams\char"76}\kern-0.3pt%
$\raise1pt\hbox{\lams\char"76}$}}\hfil}}

\def\title#1{\par
{\baselineskip1.5\baselineskip\rightskip0pt plus 5truecm
\leavevmode\vskip0truecm\noindent\font\BF=cmbx10 scaled \magstep2\BF #1\par}
\vskip1truecm
\leftline{\font\CSC=cmcsc10
{\CSC Friedrich Knop}}
\leftline{Department of Mathematics, Rutgers University, New Brunswick NJ
08903, USA}
\leftline{knop@math.rutgers.edu}
\vskip1truecm
\par}

\def\cite#1{\expandafter\ifx\csname#1\endcsname\relax
{\bf?}\immediate\write16{#1 ist nicht definiert!}\else\csname#1\endcsname\fi}
\def\expandwrite#1#2{\edef\next{\write#1{#2}}\next}
\def\neverexpand{\noexpand\noexpand\noexpand}
\def\strip#1\ {}
\def\ncite#1{\expandafter\ifx\csname#1\endcsname\relax
{\bf?}\immediate\write16{#1 ist nicht definiert!}\else
\expandafter\expandafter\expandafter\strip\csname#1\endcsname\fi}
\newwrite\AUX
\immediate\openout\AUX=\jobname.aux
\font\eightrm=cmr8\font\sixrm=cmr6
\font\eighti=cmmi8
\font\eightit=cmti8
\font\eightbf=cmbx8
\font\eightcsc=cmcsc10 scaled 833
\def\eightpoint{%
\textfont0=\eightrm\scriptfont0=\sixrm\def\rm{\fam0\eightrm}%
\textfont1=\eighti
\textfont\bffam=\eightbf\def\bf{\fam\bffam\eightbf}%
\textfont\itfam=\eightit\def\it{\fam\itfam\eightit}%
\def\csc{\eightcsc}%
\setbox\strutbox=\hbox{\vrule height7pt depth2pt width0pt}%
\normalbaselineskip=0,8\normalbaselineskip\normalbaselines\rm}
\newcount\absFnno\absFnno1
\write\AUX{\relFnno1}
\newif\ifMARKE\MARKEtrue
{\catcode`\@=11
\gdef\footnote{\ifMARKE\edef\@sf{\spacefactor\the\spacefactor}\/%
$^{\cite{Fn\the\absFnno}}$\@sf\fi
\MARKEtrue
\insert\footins\bgroup\eightpoint
\interlinepenalty100\let\par=\endgraf
\leftskip=0pt\rightskip=0pt
\splittopskip=10pt plus 1pt minus 1pt \floatingpenalty=20000\smallskip
\item{$^{\cite{Fn\the\absFnno}}$}%
\expandwrite\AUX{\neverexpand\ref{Fn\the\absFnno}{\neverexpand\the\relFnno}}%
\global\advance\absFnno1\write\AUX{\advance\relFnno1}%
\bgroup\strut\aftergroup\@foot\let\next}}
\skip\footins=12pt plus 2pt minus 4pt
\dimen\footins=30pc
\output={\plainoutput\immediate\write\AUX{\relFnno1}}
\newcount\Abschnitt\Abschnitt0
\def\beginsection#1. #2 \par{\advance\Abschnitt1%
\vskip0pt plus.10\vsize\penalty-250
\vskip0pt plus-.10\vsize\bigskip\vskip\parskip
\edef\TEST{\number\Abschnitt}
\expandafter\ifx\csname#1\endcsname\TEST\relax\else
\immediate\write16{#1 hat sich geaendert!}\fi
\expandwrite\AUX{\neverexpand\ref{#1}{\TEST}}
\leftline{\marginnote{#1}\bf\number\Abschnitt. \ignorespaces#2}%
\nobreak\smallskip\noindent\SATZ1\GNo0}
\def\Proof:{\par\noindent{\it Proof:}}
\def\Remark:{\ifdim\lastskip<\medskipamount\removelastskip\medskip\fi
\noindent{\bf Remark:}}
\def\Remarks:{\ifdim\lastskip<\medskipamount\removelastskip\medskip\fi
\noindent{\bf Remarks:}}
\def\Definition:{\ifdim\lastskip<\medskipamount\removelastskip\medskip\fi
\noindent{\bf Definition:}}
\def\Example:{\ifdim\lastskip<\medskipamount\removelastskip\medskip\fi
\noindent{\bf Example:}}
\def\Examples:{\ifdim\lastskip<\medskipamount\removelastskip\medskip\fi
\noindent{\bf Examples:}}
\newif\ifmarginalnotes\marginalnotesfalse
\newif\ifmarginalwarnings\marginalwarningstrue

\def\marginnote#1{\ifmarginalnotes\hbox to 0pt{\eightpoint\hss #1\ }\fi}

\def\strutdepth{\dp\strutbox}
\def\Randbem#1#2{\ifmarginalwarnings
{#1}\strut
\setbox0=\vtop{\eightpoint
\rightskip=0pt plus 6mm\hfuzz=3pt\hsize=16mm\noindent\leavevmode#2}%
\vadjust{\kern-\strutdepth
\vtop to \strutdepth{\kern-\ht0
\hbox to \hsize{\kern-16mm\kern-6pt\box0\kern6pt\hfill}\vss}}\fi}

\def\Zitat!{\Randbem{\bf?}{\bf Zitat}}

\newcount\SATZ\SATZ1
\def\proclaim #1. #2\par{\ifdim\lastskip<\medskipamount\removelastskip
\medskip\fi
\noindent{\bf#1.\ }{\it#2}\Par
\ifdim\lastskip<\medskipamount\removelastskip\goodbreak\medskip\fi}
\def\Aussage#1{\expandafter\def\csname#1\endcsname##1.{\resetitem
\ifx?##1?\relax\else
\edef\TEST{#1\penalty10000\ \number\Abschnitt.\number\SATZ}
\expandafter\ifx\csname##1\endcsname\TEST\relax\else
\immediate\write16{##1 hat sich geaendert!}\fi
\expandwrite\AUX{\neverexpand\ref{##1}{\TEST}}\fi
\proclaim {\marginnote{##1}\number\Abschnitt.\number\SATZ. #1\global\advance\SATZ1}.}}
\Aussage{Theorem}
\Aussage{Proposition}
\Aussage{Corollary}
\Aussage{Lemma}
\font\la=lasy10
\def\strich{\hbox{$\vcenter{\hbox
to 1pt{\leaders\hrule height -0,2pt depth 0,6pt\hfil}}$}}
\def\dashedrightarrow{\hbox{%
\hbox to 0,5cm{\leaders\hbox to 2pt{\hfil\strich\hfil}\hfil}%
\kern-2pt\hbox{\la\char\string"29}}}

\def\Bindestrich{\penalty10000-\hskip0pt}
\let\_=\Bindestrich
\def\.{{\sfcode`.=1000.}}

\def\Par{\par}
\def\:={\mathrel{\raise0,9pt\hbox{.}\kern-2,77779pt
\raise3pt\hbox{.}\kern-2,5pt=}}
\def\=:{\mathrel{=\kern-2,5pt\raise0,9pt\hbox{.}\kern-2,77779pt
\raise3pt\hbox{.}}} \def\mod{/\mskip-5mu/}
\def\into{\hookrightarrow}
\def\pfeil{\rightarrow}
\def\untenPf{\downarrow}

\def\pf#1{\buildrel#1\over\rightarrow}

\def\Ugleich{\hbox{$\cup$\kern.5pt\vrule depth -0.5pt}}
\def\|#1|{\mathop{\rm#1}\nolimits}
\def\<{\langle}
\def\>{\rangle}
\let\Times=\times
\def\times{\mathop{\Times}}
\let\Otimes=\otimes
\def\otimes{\mathop{\Otimes}}
\catcode`\@=11
\def\hex#1{\ifcase#1 0\or1\or2\or3\or4\or5\or6\or7\or8\or9\or A\or B\or
C\or D\or E\or F\else\message{Warnung: Setze hex#1=0}0\fi}
\def\fontdef#1:#2,#3,#4.{%
\alloc@8\fam\chardef\sixt@@n\FAM
\ifx!#2!\else\expandafter\font\csname text#1\endcsname=#2
\textfont\the\FAM=\csname text#1\endcsname\fi
\ifx!#3!\else\expandafter\font\csname script#1\endcsname=#3
\scriptfont\the\FAM=\csname script#1\endcsname\fi
\ifx!#4!\else\expandafter\font\csname scriptscript#1\endcsname=#4
\scriptscriptfont\the\FAM=\csname scriptscript#1\endcsname\fi
\expandafter\edef\csname #1\endcsname{\fam\the\FAM\csname text#1\endcsname}
\expandafter\edef\csname hex#1fam\endcsname{\hex\FAM}}
\catcode`\@=12 

\fontdef Ss:cmss10,,.
\fontdef Fr:eufm10,eufm7,eufm5.
\def\fa{{\Fr a}}

\def\fc{{\Fr c}}
\def\fd{{\Fr d}}

\def\fg{{\Fr g}}
\def\fh{{\Fr h}}

\def\fl{{\Fr l}}

\def\fo{{\Fr o}}
\def\fp{{\Fr p}}

\def\fs{{\Fr s}}
\def\ft{{\Fr t}}

\fontdef bbb:msbm10,msbm7,msbm5.
\fontdef mbf:cmmib10,cmmib7,.
\fontdef msa:msam10,msam7,msam5.
\def\CC{{\bbb C}}

\def\cD{{\cal D}}

\def\cO{{\cal O}}\def\cP{{\cal P}}

\def\cU{{\cal U}}\def\cW{{\cal W}}
\def\cZ{{\cal Z}}
\mathchardef\leer=\string"0\hexbbbfam3F
\mathchardef\subsetneq=\string"3\hexbbbfam24
\mathchardef\semidir=\string"2\hexbbbfam6E
\mathchardef\dirsemi=\string"2\hexbbbfam6F
\mathchardef\haken=\string"2\hexmsafam78
\mathchardef\auf=\string"3\hexmsafam10
\let\OL=\overline
\def\overline#1{{\hskip1pt\OL{\hskip-1pt#1\hskip-.3pt}\hskip.3pt}}


\def\Uq{{\overline{U}}}

%
\newdimen\Parindent
\Parindent=\parindent

\def\textindent#1{\noindent\hskip\Parindent\llap{#1\enspace }\ignorespaces}
\def\itemitem{\par\indent\hangindent2\Parindent\textindent}


\abovedisplayskip 9.0pt plus 3.0pt minus 3.0pt
\belowdisplayskip 9.0pt plus 3.0pt minus 3.0pt
\newdimen\Grenze\Grenze2\Parindent\advance\Grenze1em
\newdimen\Breite
\newbox\DpBox
\def\NewDisplay#1
#2$${\Breite\hsize\advance\Breite-\hangindent
\setbox\DpBox=\hbox{\hskip2\Parindent$\displaystyle{%
\ifx0#1\relax\else\eqno{#1}\fi#2}$}%
\ifnum\predisplaysize<\Grenze\abovedisplayskip\abovedisplayshortskip
\belowdisplayskip\belowdisplayshortskip\fi
\global\futurelet\nexttok\WEITER}
\def\WEITER{\ifx\nexttok\qed\expandafter\leftQEDdisplay
\else\leftdisplay\fi}
\def\leftdisplay{\hskip-\hangindent\leftline{\box\DpBox}$$}
\def\leftQEDdisplay{\hskip-\hangindent
\line{\copy\DpBox\hfill\lower\dp\DpBox\copy\QEDbox}%
\belowdisplayskip0pt$$\bigskip\let\nexttok=}
\everydisplay{\NewDisplay}
\newcount\GNo\GNo=0
\newcount\maxEqNo\maxEqNo=0
\def\eqno#1{%
\global\advance\GNo1
\edef\FTEST{$(\number\Abschnitt.\number\GNo)$}
\ifx?#1?\relax\else
\ifnum#1>\maxEqNo\global\maxEqNo=#1\fi%
\expandafter\ifx\csname E#1\endcsname\FTEST\relax\else
\immediate\write16{E#1 hat sich geaendert!}\fi
\expandwrite\AUX{\neverexpand\ref{E#1}{\FTEST}}\fi
\llap{\hbox to 40pt{\marginnote{#1}\FTEST\hfill}}}

\catcode`@=11
\def\eqalignno#1{\null\!\!\vcenter{\openup\jot\m@th\ialign{\eqno{##}\hfil
&\strut\hfil$\displaystyle{##}$&$\displaystyle{{}##}$\hfil\crcr#1\crcr}}\,}
\catcode`@=12

\newbox\QEDbox
\newbox\nichts\setbox\nichts=\vbox{}\wd\nichts=2mm\ht\nichts=2mm
\setbox\QEDbox=\hbox{\vrule\vbox{\hrule\copy\nichts\hrule}\vrule}
\def\qed{\leavevmode\unskip\hfil\null\nobreak\hfill\copy\QEDbox\medbreak}
\newdimen\HIindent
\newbox\HIbox
\def\setHI#1{\setbox\HIbox=\hbox{#1}\HIindent=\wd\HIbox}
\def\HI#1{\par\hangindent\HIindent\hangafter=0\noindent\leavevmode
\llap{\hbox to\HIindent{#1\hfil}}\ignorespaces}

\newdimen\maxSpalbr
\newdimen\altSpalbr
\newcount\Zaehler


\newif\ifxxx

{\catcode`/=\active

\gdef\beginrefs{%
\xxxfalse
\catcode`/=\active
\def/{\string/\ifxxx\hskip0pt\fi}
\def\TText##1{{\xxxtrue\tt##1}}
\expandafter\ifx\csname Spaltenbreite\endcsname\relax
\def\Spaltenbreite{1cm}\immediate\write16{Spaltenbreite undefiniert!}\fi
\expandafter\altSpalbr\Spaltenbreite
\maxSpalbr0pt
\gdef\alt{}
\def\\##1\relax{%
\gdef\neu{##1}\ifx\alt\neu\global\advance\Zaehler1\else
\xdef\alt{\neu}\global\Zaehler=1\fi\xdef\SigText{##1\the\Zaehler}}
\def\L|Abk:##1|Sig:##2|Au:##3|Tit:##4|Zs:##5|Bd:##6|S:##7|J:##8|xxx:##9||{%
\def\SigText{##2}\global\setbox0=\hbox{##2\relax}
\edef\TEST{[\SigText]}
\expandafter\ifx\csname##1\endcsname\TEST\relax\else
\immediate\write16{##1 hat sich geaendert!}\fi
\expandwrite\AUX{\neverexpand\ref{##1}{\TEST}}
\setHI{[\SigText]\ }
\ifnum\HIindent>\maxSpalbr\maxSpalbr\HIindent\fi
\ifnum\HIindent<\altSpalbr\HIindent\altSpalbr\fi
\HI{\marginnote{##1}[\SigText]}
\ifx-##3\relax\else{##3}: \fi
\ifx-##4\relax\else{##4}{\sfcode`.=3000.} \fi
\ifx-##5\relax\else{\it ##5\/} \fi
\ifx-##6\relax\else{\bf ##6} \fi
\ifx-##8\relax\else({##8})\fi
\ifx-##7\relax\else, {##7}\fi
\ifx-##9\relax\else, \TText{##9}\fi\Par}
\def\B|Abk:##1|Sig:##2|Au:##3|Tit:##4|Reihe:##5|Verlag:##6|Ort:##7|J:##8|xxx:##9||{%
\def\SigText{##2}\global\setbox0=\hbox{##2\relax}
\edef\TEST{[\SigText]}
\expandafter\ifx\csname##1\endcsname\TEST\relax\else
\immediate\write16{##1 hat sich geaendert!}\fi
\expandwrite\AUX{\neverexpand\ref{##1}{\TEST}}
\setHI{[\SigText]\ }
\ifnum\HIindent>\maxSpalbr\maxSpalbr\HIindent\fi
\ifnum\HIindent<\altSpalbr\HIindent\altSpalbr\fi
\HI{\marginnote{##1}[\SigText]}
\ifx-##3\relax\else{##3}: \fi
\ifx-##4\relax\else{##4}{\sfcode`.=3000.} \fi
\ifx-##5\relax\else{(##5)} \fi
\ifx-##7\relax\else{##7:} \fi
\ifx-##6\relax\else{##6}\fi
\ifx-##8\relax\else{ ##8}\fi
\ifx-##9\relax\else, \TText{##9}\fi\Par}
\def\Pr|Abk:##1|Sig:##2|Au:##3|Artikel:##4|Titel:##5|Hgr:##6|Reihe:{%
\def\SigText{##2}\global\setbox0=\hbox{##2\relax}
\edef\TEST{[\SigText]}
\expandafter\ifx\csname##1\endcsname\TEST\relax\else
\immediate\write16{##1 hat sich geaendert!}\fi
\expandwrite\AUX{\neverexpand\ref{##1}{\TEST}}
\setHI{[\SigText]\ }
\ifnum\HIindent>\maxSpalbr\maxSpalbr\HIindent\fi
\ifnum\HIindent<\altSpalbr\HIindent\altSpalbr\fi
\HI{\marginnote{##1}[\SigText]}
\ifx-##3\relax\else{##3}: \fi
\ifx-##4\relax\else{##4}{\sfcode`.=3000.} \fi
\ifx-##5\relax\else{In: \it ##5}. \fi
\ifx-##6\relax\else{(##6)} \fi\PrII}
\def\PrII##1|Bd:##2|Verlag:##3|Ort:##4|S:##5|J:##6|xxx:##7||{%
\ifx-##1\relax\else{##1} \fi
\ifx-##2\relax\else{\bf ##2}, \fi
\ifx-##4\relax\else{##4:} \fi
\ifx-##3\relax\else{##3} \fi
\ifx-##6\relax\else{##6}\fi
\ifx-##5\relax\else{, ##5}\fi
\ifx-##7\relax\else, \TText{##7}\fi\Par}
\bgroup
\baselineskip12pt
\parskip2.5pt plus 1pt
\hyphenation{Hei-del-berg Sprin-ger}
\sfcode`.=1000
\beginsection References. References

}}

\def\endrefs{%
\expandwrite\AUX{\neverexpand\ref{Spaltenbreite}{\the\maxSpalbr}}
\ifnum\maxSpalbr=\altSpalbr\relax\else
\immediate\write16{Spaltenbreite hat sich geaendert!}\fi
\egroup\write16{Letzte Gleichung: E\the\maxEqNo}
\write16{Letzte Aufzaehlung: I\the\maxItemcount}}







\def\setRevDate $#1 #2 #3${#2}
\def\TeXdrawId{\setRevDate $Date: 1995/12/19 16:40:42 $ TeXdraw V2R0}
\chardef\catamp=\the\catcode`\@
\catcode`\@=11
\long
\def\centertexdraw #1{\hbox to \hsize{\hss
\btexdraw #1\etexdraw
\hss}}
\def\btexdraw {\x@pix=0             \y@pix=0
\x@segoffpix=\x@pix  \y@segoffpix=\y@pix
\t@exdrawdef
\setbox\t@xdbox=\vbox\bgroup\offinterlineskip
\global\d@bs=0
\global\t@extonlytrue
\global\p@osinitfalse
\s@avemove \x@pix \y@pix
\m@pendingfalse
\global\p@osinitfalse
\p@athfalse
\the\everytexdraw}
\def\etexdraw {\ift@extonly \else
\t@drclose
\fi
\egroup
\ifdim \wd\t@xdbox>0pt
\t@xderror {TeXdraw box non-zero size,
possible extraneous text}%
\fi
\vbox {\offinterlineskip
\pixtobp \xminpix \l@lxbp  \pixtobp \yminpix \l@lybp
\pixtobp \xmaxpix \u@rxbp  \pixtobp \ymaxpix \u@rybp
\hbox{\t@xdinclude 
[{\l@lxbp},{\l@lybp}][{\u@rxbp},{\u@rybp}]{\p@sfile}}%
\pixtodim \xminpix \t@xpos  \pixtodim \yminpix \t@ypos
\kern \t@ypos
\hbox {\kern -\t@xpos
\box\t@xdbox
\kern \t@xpos}%
\kern -\t@ypos\relax}}
\def\drawdim #1 {\def\d@dim{#1\relax}}
\def\setunitscale #1 {\edef\u@nitsc{#1}%
\realmult \u@nitsc \s@egsc \d@sc}
\def\relunitscale #1 {\realmult {#1}\u@nitsc \u@nitsc
\realmult \u@nitsc \s@egsc \d@sc}
\def\setsegscale #1 {\edef\s@egsc {#1}%
\realmult \u@nitsc \s@egsc \d@sc}
\def\relsegscale #1 {\realmult {#1}\s@egsc \s@egsc
\realmult \u@nitsc \s@egsc \d@sc}
\def\bsegment {\ifp@ath
\f@lushbs
\f@lushmove
\fi
\begingroup
\x@segoffpix=\x@pix
\y@segoffpix=\y@pix
\setsegscale 1
\global\advance \d@bs by 1\relax}
\def\esegment {\endgroup
\ifnum \d@bs=0
\writetx {es}%
\else
\global\advance \d@bs by -1
\fi}
\def\savecurrpos (#1 #2){\getsympos (#1 #2)\a@rgx\a@rgy
\s@etcsn \a@rgx {\the\x@pix}%
\s@etcsn \a@rgy {\the\y@pix}}
\def\savepos (#1 #2)(#3 #4){\getpos (#1 #2)\a@rgx\a@rgy
\coordtopix \a@rgx \t@pixa
\advance \t@pixa by \x@segoffpix
\coordtopix \a@rgy \t@pixb
\advance \t@pixb by \y@segoffpix
\getsympos (#3 #4)\a@rgx\a@rgy
\s@etcsn \a@rgx {\the\t@pixa}%
\s@etcsn \a@rgy {\the\t@pixb}}
\def\linewd #1 {\coordtopix {#1}\t@pixa
\f@lushbs
\writetx {\the\t@pixa\space sl}}
\def\setgray #1 {\f@lushbs
\writetx {#1 sg}}
\def\lpatt (#1){\listtopix (#1)\p@ixlist
\f@lushbs
\writetx {[\p@ixlist] sd}}
\def\lvec (#1 #2){\getpos (#1 #2)\a@rgx\a@rgy
\s@etpospix \a@rgx \a@rgy
\writeps {\the\x@pix\space \the\y@pix\space lv}}
\def\rlvec (#1 #2){\getpos (#1 #2)\a@rgx\a@rgy
\r@elpospix \a@rgx \a@rgy
\writeps {\the\x@pix\space \the\y@pix\space lv}}
\def\move (#1 #2){\getpos (#1 #2)\a@rgx\a@rgy
\s@etpospix \a@rgx \a@rgy
\s@avemove \x@pix \y@pix}
\def\rmove (#1 #2){\getpos (#1 #2)\a@rgx\a@rgy
\r@elpospix \a@rgx \a@rgy
\s@avemove \x@pix \y@pix}
\def\lcir r:#1 {\coordtopix {#1}\t@pixa
\writeps {\the\t@pixa\space cr}%
\r@elupd \t@pixa \t@pixa
\r@elupd {-\t@pixa}{-\t@pixa}}
\def\fcir f:#1 r:#2 {\coordtopix {#2}\t@pixa
\writeps {\the\t@pixa\space #1 fc}%
\r@elupd \t@pixa \t@pixa
\r@elupd {-\t@pixa}{-\t@pixa}}
\def\lellip rx:#1 ry:#2 {\coordtopix {#1}\t@pixa
\coordtopix {#2}\t@pixb
\writeps {\the\t@pixa\space \the\t@pixb\space el}%
\r@elupd \t@pixa \t@pixb
\r@elupd {-\t@pixa}{-\t@pixb}}
\def\fellip f:#1 rx:#2 ry:#3 {\coordtopix {#2}\t@pixa
\coordtopix {#3}\t@pixb
\writeps {\the\t@pixa\space \the\t@pixb\space #1 fe}%
\r@elupd \t@pixa \t@pixb
\r@elupd {-\t@pixa}{-\t@pixb}}
\def\larc r:#1 sd:#2 ed:#3 {\coordtopix {#1}\t@pixa
\writeps {\the\t@pixa\space #2 #3 ar}}
\def\ifill f:#1 {\writeps {#1 fl}}
\def\lfill f:#1 {\writeps {#1 fp}}
\def\htext #1{\def\testit {#1}%
\ifx \testit\l@paren
\let\next=\h@move
\else
\let\next=\h@text
\fi
\next {#1}}
\def\rtext td:#1 #2{\def\testit {#2}%
\ifx \testit\l@paren
\let\next=\r@move
\else
\let\next=\r@text
\fi
\next td:#1 {#2}}

\def\textref h:#1 v:#2 {\ifx #1R%
\edef\l@stuff {\hss}\edef\r@stuff {}%
\else
\ifx #1C%
\edef\l@stuff {\hss}\edef\r@stuff {\hss}%
\else
\edef\l@stuff {}\edef\r@stuff {\hss}%
\fi
\fi
\ifx #2T%
\edef\t@stuff {}\edef\b@stuff {\vss}%
\else
\ifx #2C%
\edef\t@stuff {\vss}\edef\b@stuff {\vss}%
\else
\edef\t@stuff {\vss}\edef\b@stuff {}%
\fi
\fi}
\def\avec (#1 #2){\getpos (#1 #2)\a@rgx\a@rgy
\s@etpospix \a@rgx \a@rgy
\writeps {\the\x@pix\space \the\y@pix\space (\a@type)
\the\a@lenpix\space \the\a@widpix\space av}}
\def\ravec (#1 #2){\getpos (#1 #2)\a@rgx\a@rgy
\r@elpospix \a@rgx \a@rgy
\writeps {\the\x@pix\space \the\y@pix\space (\a@type)
\the\a@lenpix\space \the\a@widpix\space av}}
\def\arrowheadsize l:#1 w:#2 {\coordtopix{#1}\a@lenpix
\coordtopix{#2}\a@widpix}
\def\arrowheadtype t:#1 {\edef\a@type{#1}}
\def\clvec (#1 #2)(#3 #4)(#5 #6)%
{\getpos (#1 #2)\a@rgx\a@rgy
\coordtopix \a@rgx\t@pixa
\advance \t@pixa by \x@segoffpix
\coordtopix \a@rgy\t@pixb
\advance \t@pixb by \y@segoffpix
\getpos (#3 #4)\a@rgx\a@rgy
\coordtopix \a@rgx\t@pixc
\advance \t@pixc by \x@segoffpix
\coordtopix \a@rgy\t@pixd
\advance \t@pixd by \y@segoffpix
\getpos (#5 #6)\a@rgx\a@rgy
\s@etpospix \a@rgx \a@rgy
\writeps {\the\t@pixa\space \the\t@pixb\space
\the\t@pixc\space \the\t@pixd\space
\the\x@pix\space \the\y@pix\space cv}}
\def\drawbb {\bsegment
\drawdim bp
\linewd 0.24
\setunitscale {\p@sfactor}
\writeps {\the\xminpix\space \the\yminpix\space mv}%
\writeps {\the\xminpix\space \the\ymaxpix\space lv}%
\writeps {\the\xmaxpix\space \the\ymaxpix\space lv}%
\writeps {\the\xmaxpix\space \the\yminpix\space lv}%
\writeps {\the\xminpix\space \the\yminpix\space lv}%
\esegment}
\def\getpos (#1 #2)#3#4{\g@etargxy #1 #2 {} \\#3#4%
\c@heckast #3%
\ifa@st
\g@etsympix #3\t@pixa
\advance \t@pixa by -\x@segoffpix
\pixtocoord \t@pixa #3%
\fi
\c@heckast #4%
\ifa@st
\g@etsympix #4\t@pixa
\advance \t@pixa by -\y@segoffpix
\pixtocoord \t@pixa #4%
\fi}
\def\getsympos (#1 #2)#3#4{\g@etargxy #1 #2 {} \\#3#4%
\c@heckast #3%
\ifa@st \else
\t@xderror {TeXdraw: invalid symbolic coordinate}%
\fi
\c@heckast #4%
\ifa@st \else
\t@xderror {TeXdraw: invalid symbolic coordinate}%
\fi}
\def\listtopix (#1)#2{\def #2{}%
\edef\l@ist {#1 }%
\m@oretrue
\loop
\expandafter\g@etitem \l@ist \\\a@rgx\l@ist
\a@pppix \a@rgx #2%
\ifx \l@ist\empty
\m@orefalse
\fi
\ifm@ore
\repeat}
\def\realmult #1#2#3{\dimen0=#1pt
\dimen2=#2\dimen0
\edef #3{\expandafter\c@lean\the\dimen2}}
\def\intdiv #1#2#3{\t@counta=#1
\t@countb=#2
\ifnum \t@countb<0
\t@counta=-\t@counta
\t@countb=-\t@countb
\fi
\t@countd=1
\ifnum \t@counta<0
\t@counta=-\t@counta
\t@countd=-1
\fi
\t@countc=\t@counta  \divide \t@countc by \t@countb
\t@counte=\t@countc  \multiply \t@counte by \t@countb
\advance \t@counta by -\t@counte
\t@counte=-1
\loop
\advance \t@counte by 1
\ifnum \t@counte<16
\multiply \t@countc by 2
\multiply \t@counta by 2
\ifnum \t@counta<\t@countb \else
\advance \t@countc by 1
\advance \t@counta by -\t@countb
\fi
\repeat
\divide \t@countb by 2
\ifnum \t@counta<\t@countb
\advance \t@countc by 1
\fi
\ifnum \t@countd<0
\t@countc=-\t@countc
\fi
\dimen0=\t@countc sp
\edef #3{\expandafter\c@lean\the\dimen0}}
\def\coordtopix #1#2{\dimen0=#1\d@dim
\dimen2=\d@sc\dimen0
\t@counta=\dimen2
\t@countb=\s@ppix
\divide \t@countb by 2
\ifnum \t@counta<0
\advance \t@counta by -\t@countb
\else
\advance \t@counta by \t@countb
\fi
\divide \t@counta by \s@ppix
#2=\t@counta}
\def\pixtocoord #1#2{\t@counta=#1%
\multiply \t@counta by \s@ppix
\dimen0=\d@sc\d@dim
\t@countb=\dimen0
\intdiv \t@counta \t@countb #2}
\def\pixtodim #1#2{\t@countb=#1%
\multiply \t@countb by \s@ppix
#2=\t@countb sp\relax}
\def\pixtobp #1#2{\dimen0=\p@sfactor pt
\t@counta=\dimen0
\multiply \t@counta by #1%
\ifnum \t@counta < 0
\advance \t@counta by -32768
\else
\advance \t@counta by 32768
\fi
\divide \t@counta by 65536
\edef #2{\the\t@counta}}
\newcount\t@counta    \newcount\t@countb
\newcount\t@countc    \newcount\t@countd
\newcount\t@counte
\newcount\t@pixa      \newcount\t@pixb
\newcount\t@pixc      \newcount\t@pixd
\newdimen\t@xpos      \newdimen\t@ypos
\newcount\xminpix      \newcount\xmaxpix
\newcount\yminpix      \newcount\ymaxpix
\newcount\a@lenpix     \newcount\a@widpix
\newcount\x@pix        \newcount\y@pix
\newcount\x@segoffpix  \newcount\y@segoffpix
\newcount\x@savepix    \newcount\y@savepix
\newcount\s@ppix
\newcount\d@bs
\newcount\t@xdnum
\global\t@xdnum=0
\newbox\t@xdbox
\newwrite\drawfile
\newif\ifm@pending
\newif\ifp@ath
\newif\ifa@st
\newif\ifm@ore
\newif \ift@extonly
\newif\ifp@osinit
\newtoks\everytexdraw
\def\l@paren{(}
\def\a@st{*}
\catcode`\%=12
\def\p@b {
\catcode`\%=14
\catcode`\{=12  \catcode`\}=12  \catcode`\u=1 \catcode`\v=2
\def\l@br u{v  \def\r@br u}v
\catcode `\{=1  \catcode`\}=2   \catcode`\u=11 \catcode`\v=11
{\catcode`\p=12 \catcode`\t=12
\gdef\c@lean #1pt{#1}}
\def\sppix#1/#2 {\dimen0=1#2 \s@ppix=\dimen0
\t@counta=#1%
\divide \t@counta by 2
\advance \s@ppix by \t@counta
\divide \s@ppix by #1%
\t@counta=\s@ppix
\multiply \t@counta by 65536
\advance \t@counta by 32891
\divide \t@counta by 65782
\dimen0=\t@counta sp
\edef\p@sfactor {\expandafter\c@lean\the\dimen0}}
\def\g@etargxy #1 #2 #3 #4\\#5#6{\def #5{#1}%
\ifx #5\empty
\g@etargxy #2 #3 #4 \\#5#6
\else
\def #6{#2}%
\def\next {#3}%
\ifx \next\empty \else
\t@xderror {TeXdraw: invalid coordinate}%
\fi
\fi}
\def\c@heckast #1{\expandafter
\c@heckastll #1\\}
\def\c@heckastll #1#2\\{\def\testit {#1}%
\ifx \testit\a@st
\a@sttrue
\else
\a@stfalse
\fi}
\def\g@etsympix #1#2{\expandafter
\ifx \csname #1\endcsname \relax
\t@xderror {TeXdraw: undefined symbolic coordinate}%
\fi
#2=\csname #1\endcsname}
\def\s@etcsn #1#2{\expandafter
\xdef\csname#1\endcsname {#2}}
\def\g@etitem #1 #2\\#3#4{\edef #4{#2}\edef #3{#1}}
\def\a@pppix #1#2{\edef\next {#1}%
\ifx \next\empty \else
\coordtopix {#1}\t@pixa
\ifx #2\empty
\edef #2{\the\t@pixa}%
\else
\edef #2{#2 \the\t@pixa}%
\fi
\fi}
\def\s@etpospix #1#2{\coordtopix {#1}\x@pix
\advance \x@pix by \x@segoffpix
\coordtopix {#2}\y@pix
\advance \y@pix by \y@segoffpix
\u@pdateminmax \x@pix \y@pix}
\def\r@elpospix #1#2{\coordtopix {#1}\t@pixa
\advance \x@pix by \t@pixa
\coordtopix {#2}\t@pixa
\advance \y@pix by \t@pixa
\u@pdateminmax \x@pix \y@pix}
\def\r@elupd #1#2{\t@counta=\x@pix
\advance\t@counta by #1%
\t@countb=\y@pix
\advance\t@countb by #2%
\u@pdateminmax \t@counta \t@countb}
\def\u@pdateminmax #1#2{\ifnum #1>\xmaxpix
\global\xmaxpix=#1%
\fi
\ifnum #1<\xminpix
\global\xminpix=#1%
\fi
\ifnum #2>\ymaxpix
\global\ymaxpix=#2%
\fi
\ifnum #2<\yminpix
\global\yminpix=#2%
\fi}
\def\s@avemove #1#2{\x@savepix=#1\y@savepix=#2%
\m@pendingtrue
\ifp@osinit \else
\global\p@osinittrue
\global\xminpix=\x@savepix \global\yminpix=\y@savepix
\global\xmaxpix=\x@savepix \global\ymaxpix=\y@savepix
\fi}
\def\f@lushmove {\global\p@osinittrue
\ifm@pending
\writetx {\the\x@savepix\space \the\y@savepix\space mv}%
\m@pendingfalse
\p@athfalse
\fi}
\def\f@lushbs {\loop
\ifnum \d@bs>0
\writetx {bs}%
\global\advance \d@bs by -1
\repeat}
\def\h@move #1#2 #3)#4{\move (#2 #3)%
\h@text {#4}}
\def\h@text #1{\pixtodim \x@pix \t@xpos
\pixtodim \y@pix \t@ypos
\vbox to 0pt{\normalbaselines
\t@stuff
\kern -\t@ypos
\hbox to 0pt{\l@stuff
\kern \t@xpos
\hbox {#1}%
\kern -\t@xpos
\r@stuff}%
\kern \t@ypos
\b@stuff\relax}}
\def\r@move td:#1 #2#3 #4)#5{\move (#3 #4)%
\r@text td:#1 {#5}}
\def\r@text td:#1 #2{\vbox to 0pt{\pixtodim \x@pix \t@xpos
\pixtodim \y@pix \t@ypos
\kern -\t@ypos
\hbox to 0pt{\kern \t@xpos
\rottxt {#1}{\z@sb {#2}}%
\hss}%
\vss}}
\def\z@sb #1{\vbox to 0pt{\normalbaselines
\t@stuff
\hbox to 0pt{\l@stuff \hbox {#1}\r@stuff}%
\b@stuff}}
\ifx \rotatebox\@undefined
\def\rottxt #1#2{\bgroup
#2%
\egroup}
\else
\let\rottxt=\rotatebox
\fi
\ifx \t@xderror\@undefined
\let\t@xderror=\errmessage
\fi
\def\t@exdrawdef {\sppix 300/in
\drawdim in
\edef\u@nitsc {1}%
\setsegscale 1
\arrowheadsize l:0.16 w:0.08
\arrowheadtype t:T
\textref h:L v:B }
\ifx \includegraphics\@undefined
\def\t@xdinclude [#1,#2][#3,#4]#5{%
\begingroup
\message {<#5>}%
\leavevmode
\t@counta=-#1%
\t@countb=-#2%
\setbox0=\hbox{%
\includegraphics{#5}}%
\t@ypos=#4 bp%
\advance \t@ypos by -#2 bp%
\t@xpos=#3 bp%
\advance \t@xpos by -#1 bp%
\dp0=0pt \ht0=\t@ypos  \wd0=\t@xpos
\box0%
\endgroup}
\else
\let\t@xdinclude=\includegraphics
\fi
\def\writeps #1{\f@lushbs
\f@lushmove
\p@athtrue
\writetx {#1}}
\def\writetx #1{\ift@extonly
\global\t@extonlyfalse
\t@xdpsfn \p@sfile
\t@dropen \p@sfile
\fi
\w@rps {#1}}
\def\w@rps #1{\immediate\write\drawfile {#1}}
\def\t@xdpsfn #1{%
\global\advance \t@xdnum by 1
\ifnum \t@xdnum<10
\xdef #1{\jobname.ps\the\t@xdnum}
\else
\xdef #1{\jobname.p\the\t@xdnum}
\fi
}
\def\t@dropen #1{%
\immediate\openout\drawfile=#1%
\w@rps {\p@b PS-Adobe-3.0 EPSF-3.0}%
\w@rps {\p@p BoundingBox: (atend)}%
\w@rps {\p@p Title: TeXdraw drawing: #1}%
\w@rps {\p@p Pages: 1}%
\w@rps {\p@p Creator: \TeXdrawId}%
\w@rps {\p@p CreationDate: \the\year/\the\month/\the\day}%
\w@rps {50 dict begin}%
\w@rps {/mv {stroke moveto} def}%
\w@rps {/lv {lineto} def}%
\w@rps {/st {currentpoint stroke moveto} def}%
\w@rps {/sl {st setlinewidth} def}%
\w@rps {/sd {st 0 setdash} def}%
\w@rps {/sg {st setgray} def}%
\w@rps {/bs {gsave} def /es {stroke grestore} def}%
\w@rps {/fl \l@br gsave setgray fill grestore}%
\w@rps    { currentpoint newpath moveto\r@br\space def}%
\w@rps {/fp {gsave setgray fill grestore st} def}%
\w@rps {/cv {curveto} def}%
\w@rps {/cr \l@br gsave currentpoint newpath 3 -1 roll 0 360 arc}%
\w@rps    { stroke grestore\r@br\space def}%
\w@rps {/fc \l@br gsave setgray currentpoint newpath}%
\w@rps    { 3 -1 roll 0 360 arc fill grestore\r@br\space def}%
\w@rps {/ar {gsave currentpoint newpath 5 2 roll arc stroke grestore} def}%
\w@rps {/el \l@br gsave /svm matrix currentmatrix def}%
\w@rps    { currentpoint translate scale newpath 0 0 1 0 360 arc}%
\w@rps    { svm setmatrix stroke grestore\r@br\space def}%
\w@rps {/fe \l@br gsave setgray currentpoint translate scale newpath}%
\w@rps    { 0 0 1 0 360 arc fill grestore\r@br\space def}%
\w@rps {/av \l@br /hhwid exch 2 div def /hlen exch def}%
\w@rps    { /ah exch def /tipy exch def /tipx exch def}%
\w@rps    { currentpoint /taily exch def /tailx exch def}%
\w@rps    { /dx tipx tailx sub def /dy tipy taily sub def}%
\w@rps    { /alen dx dx mul dy dy mul add sqrt def}%
\w@rps    { /blen alen hlen sub def}%
\w@rps    { gsave tailx taily translate dy dx atan rotate}%
\w@rps    { (V) ah ne {blen 0 gt {blen 0 lineto} if} {alen 0 lineto} ifelse}%
\w@rps    { stroke blen hhwid neg moveto alen 0 lineto blen hhwid lineto}%
\w@rps    { (T) ah eq {closepath} if}%
\w@rps    { (W) ah eq {gsave 1 setgray fill grestore closepath} if}%
\w@rps    { (F) ah eq {fill} {stroke} ifelse}%
\w@rps    { grestore tipx tipy moveto\r@br\space def}%
\w@rps {\p@sfactor\space \p@sfactor\space scale}%
\w@rps {1 setlinecap 1 setlinejoin}%
\w@rps {3 setlinewidth [] 0 setdash}%
\w@rps {0 0 moveto}%
}
\def\t@drclose {%
\bgroup
\w@rps {stroke end showpage}%
\w@rps {\p@p Trailer:}%
\pixtobp \xminpix \l@lxbp  \pixtobp \yminpix \l@lybp
\pixtobp \xmaxpix \u@rxbp  \pixtobp \ymaxpix \u@rybp
\w@rps {\p@p BoundingBox: \l@lxbp\space \l@lybp\space
\u@rxbp\space \u@rybp}%
\w@rps {\p@p EOF}%
\egroup
\immediate\closeout\drawfile
}
\catcode`\@=\catamp

\def\Rad{0.1}
\let\rrr=\Rad
\realmult2\Rad\RRad
\def\STL{0.8}
\def\factor{0.6}
\def\unit{cm}
\newdimen\Radius
\Radius=\Rad\unit
\Radius=\factor\Radius
\everytexdraw{%
\drawdim pt
\linewd 0.2
\drawdim{\unit}
\setunitscale{\factor}
}
{\catcode`\p=12\catcode`\t=12\gdef\sri#1pt{#1}}
\newbox\txt
\newdimen\Hoehe
\def\\{\cr}
\def\Zeichen#1{\htext{\BOX{#1}}}
\def\BOX#1{\setbox\txt=\vbox{\baselineskip6pt%
\halign{$\scriptstyle##$\hfil\cr#1\cr}}%
\global\edef\B{\expandafter\sri\the\wd\txt}%
\global\Hoehe=\ht\txt
\global\advance\Hoehe\dp\txt
\global\edef\H{\expandafter\sri\the\Hoehe}%
\copy\txt}
\def\punkt(#1,#2){\move(#1 #2)\fcir f:0 r:\Rad}
\def\Kreis(#1,#2){\move(#1 #2)\lcir r:\RRad}
\def\kreis(#1,#2){\move(#1 #2)\lcir r:\Rad}
\def\AstrichR(#1,#2){\move(#1 #2)\rmove(0.1 0)\rlvec(0.8 0)\rmove(0.1 0)}
\def\BstrichR(#1,#2){\move(#1 #2)\rmove(0.07071 0.07071)\rlvec(0.85858 0)
\move(#1 #2)\rmove(0.07071 -0.07071)\rlvec(0.85858 0)\rmove(0.07071 0.07071)
\rmove(-0.3 0)\rlvec(-0.2 -0.2)\rmove(0.2 0.2)\rlvec(-0.2 0.2)}
\def\CstrichR(#1,#2){\move(#1 #2)\rmove(0.07071 0.07071)\rlvec(0.85858 0)
\move(#1 #2)\rmove(0.07071 -0.07071)\rlvec(0.85858 0)\move(#1 #2)
\rmove(0.3 0)\rlvec(0.2 0.2)\rmove(-0.2 -0.2)\rlvec(0.2 -0.2)}
\def\GstrichR(#1,#2){\AstrichR(#1,#2)\move(#1 #2)\rmove(0.07071 0.07071)\rlvec(0.85858 0)
\move(#1 #2)\rmove(0.07071 -0.07071)\rlvec(0.85858 0)\move(#1 #2)
\rmove(0.3 0)\rlvec(0.2 0.2)\rmove(-0.2 -0.2)\rlvec(0.2 -0.2)}
\def\DstrichRO(#1,#2){\move(#1 #2)\rmove(0.07071 0.07071)\rlvec(0.85858 0.85858)\rmove(0.07071 0.07071)}
\def\DstrichRU(#1,#2){\move(#1 #2)\rmove(0.07071 -0.07071)\rlvec(0.85858 -0.85858)\rmove(0.07071 -0.07071)}
\def\strichRO(#1,#2){\move(#1 #2)\rlvec(1 1)}
\def\strichRU(#1,#2){\move(#1 #2)\rlvec(1 -1)}
\def\strichU(#1,#2){\move(#1 #2)\rlvec(0 -1)}
\def\labelO(#1,#2)[#3]#4{%
\punkt(#1,#2)
\textref h:C v:B
\Zeichen{#4\vrule depth 5pt width 0pt}
{\drawdim pt
\setunitscale 1
\rmove(0 {\H})
\setunitscale 0.5
\rmove({\B} 0)
\setunitscale 1
\rmove({-\B} 0)
}}
\def\KlabelU(#1,#2)[#3]#4{%
\kreis(#1,#2)
\move(#1 #2)
\rmove(#3 -0.15)
\textref h:C v:T
\Zeichen{#4}
{\drawdim pt
\setunitscale 1
\rmove(0 {-\H})
\setunitscale 0.5
\rmove({\B} 0)
\setunitscale 1
\rmove({-\B} 0)
}\ifx0#2\realmult{0.15}{\factor}\fff
\advance\Hoehe by \fff\unit
\ifdim\Hoehe > \maxT \global\maxT=\Hoehe\fi\fi}
\def\PlabelU(#1,#2)[#3]#4{%
\punkt(#1,#2)
\rmove(#3 -0.15)
\textref h:C v:T
\Zeichen{#4}
{\drawdim pt
\setunitscale 1
\rmove(0 {-\H})
\setunitscale 0.5
\rmove({\B} 0)
\setunitscale 1
\rmove({-\B} 0)
}\ifx0#2\realmult{0.15}{\factor}\fff
\advance\Hoehe by \fff\unit
\ifdim\Hoehe > \maxT \global\maxT=\Hoehe\fi\fi}
\def\PlabelR(#1,#2)[#3]#4{%
\punkt(#1,#2)
\rmove(.2 #3)
\textref h:L v:C
\Zeichen{#4}
{\drawdim pt
\setunitscale 1
\rmove({\B} 0)
\setunitscale 0.5
\rmove(0 {\H})
\setunitscale 1
\rmove(0 {-\H})
}\ifx0#2\Hoehe=0.5\Hoehe
\realmult{#3}{\factor}\fff
\advance\Hoehe by -\fff\unit
\ifdim\Hoehe > \maxT \global\maxT=\Hoehe\fi\fi}
\def\KlabelR(#1,#2)[#3]#4{%
\kreis(#1,#2)
\rmove(.2 #3)
\textref h:L v:C
\Zeichen{#4}
{\drawdim pt
\setunitscale 1
\rmove({\B} 0)
\setunitscale 0.5
\rmove(0 {\H})
\setunitscale 1
\rmove(0 {-\H})
}\ifx0#2\Hoehe=0.5\Hoehe
\realmult{#3}{\factor}\fff
\advance\Hoehe by -\fff\unit
\ifdim\Hoehe > \maxT \global\maxT=\Hoehe\fi\fi}
\def\PlabelL(#1,#2)[#3]#4{%
\punkt(#1,#2)
\rmove(-.2 #3)
\textref h:R v:C
\Zeichen{#4}
{\drawdim pt
\setunitscale 1
\rmove({-\B} 0)
\setunitscale 0.5
\rmove(0 {\H})
\setunitscale 1
\rmove(0 {-\H})
}\ifx0#2\Hoehe=0.5\Hoehe
\realmult{#3}{\factor}\fff
\advance\Hoehe by -\fff\unit
\ifdim\Hoehe > \maxT \global\maxT=\Hoehe\fi\fi}
\def\KlabelL(#1,#2)[#3]#4{%
\kreis(#1,#2)
\rmove(-.2 #3)
\textref h:R v:C
\Zeichen{#4}
{\drawdim pt
\setunitscale 1
\rmove({-\B} 0)
\setunitscale 0.5
\rmove(0 {\H})
\setunitscale 1
\rmove(0 {-\H})
}\ifx0#2\Hoehe=0.5\Hoehe
\realmult{#3}{\factor}\fff
\advance\Hoehe by -\fff\unit
\ifdim\Hoehe > \maxT \global\maxT=\Hoehe\fi\fi}
\def\Vstern(#1,#2){%
\move(#1 #2)
\rmove(-0.1 0.5)
\textref h:R v:C
\Zeichen{\textstyle*}
{\drawdim pt
\setunitscale 1
\rmove({-\B} 0)
\setunitscale 0.5
\rmove(0 {\H})
\setunitscale 1
\rmove(0 {-\H})
}}
\def\Lstern(#1,#2){%
\move(#1 #2)
\rmove(-0.8 0.4)
\textref h:C v:C
\Zeichen{\textstyle*}}
\def\Mstern(#1,#2){%
\move(#1 #2)
\rmove(0 -0.4)
\textref h:C v:C
\Zeichen{\textstyle*}}
\newdimen\maxT
\newbox\Zbox
\long\def\bdiadraw#1\ediadraw{%
\maxT=0pt\setbox\Zbox=\hbox{\btexdraw #1\etexdraw}%
\ifdim\maxT>\Radius\advance\maxT by -\Radius\fi
\lower\maxT\box\Zbox}

\def\sl{{\fs\fl}}
\def\gl{{\fg\fl}}
\def\so{{\fs\fo}}
\def\sp{{\fs\fp}}
\def\spin{{\Fr spin}}

\def\vplus{\mathop{\underline\oplus}}

\let\UL=\underline
\newdimen\uldepth

\def\underline#1{\setbox0=\hbox{$#1$}\uldepth=\dp0%
\setbox0=\hbox{$\UL{#1}$}\dp0=\uldepth\box0}

\def\usl{\underline\sl}

\fontdef sans:cmss10,,.

\def\sA#1{{\sans A}_{#1}}
\def\sB#1{{\sans B}_{#1}}
\def\sC#1{{\sans C}_{#1}}
\def\sD#1{{\sans D}_{#1}}
\def\sE#1{{\sans E}_{#1}}
\def\sG#1{{\sans G}_{#1}}

\mathchardef\ersetze=\string"3\hexmsafam20

\def\hideA#1{#1}
\def\hideWeyl#1{#1}
\def\hideEmb#1{}
\def\hideIndex#1{#1}

\newdimen\rand\rand=5pt
\def\wrapit#1{\hbox{\hskip\rand\vbox{\vskip\rand#1\vskip\rand}\hskip\rand}}
\def\boxit#1{\hbox{\vrule\vbox{\hrule\wrapit{#1}\hrule}\vrule}}

\newcount\xxx\xxx=0
\def\weiter{\global\advance\xxx1\hbox to \SpBr{$\langle\hbox{\tablename.}\number\xxx\rangle$\ \hfill}}
\def\nweiter#1&{\omit\hskip\SpBr\phantom{: }$#1$\quad\hfill&}

\newdimen\SpBr

\def\ITEM#1#2{\global\advance\itemcount1
\edef\TEXT{{\tenrm[{\bf#2}]}}%
\ifx?#1?\relax\else
\ifnum#1>\maxItemcount\global\maxItemcount=#1\fi
\expandafter\ifx\csname I#1\endcsname\TEXT\relax\else
\immediate\write16{I#1 hat sich geaendert!}\fi
\expandwrite\AUX{\neverexpand\ref{I#1}{\TEXT}}\fi
\item{\TEXT}}

\chardef\catamp=\the\catcode`\@
\catcode`\@=11
\def\realadd #1#2#3{\dimen0=#1pt
\dimen2=#2pt
\advance \dimen0 by \dimen2
\edef #3{\expandafter\c@lean\the\dimen0}}
\catcode`\@=\catamp
\let\Move=\move
\def\move(#1 #2){\realadd{#1}{.1}\Xcoord\realadd{#2}{.5}\Ycoord\Move({\Xcoord} {\Ycoord})}

\def\Nichts{\hbox{---}}


\title{Classification of multiplicity free symplectic representations}
\bigskip

{\baselineskip12pt\narrower\noindent {\bf Abstract.} Let $G$ be a
connected reductive group acting on a finite dimensional vector space
$V$. Assume that $V$ is equipped with a $G$-invariant symplectic
form. Then the ring $\cO(V)$ of polynomial functions becomes a Poisson
algebra. The ring $\cO(V)^G$ of invariants is a sub-Poisson
algebra. We call $V$ multiplicity free if $\cO(V)^G$ is Poisson
commutative, i.e., if $\{f,g\}=0$ for all invariants $f$ and
$g$. Alternatively, $G$ also acts on the Weyl algebra $\cW(V)$ and $V$
is multiplicity free if and only if the subalgebra $\cW(V)^G$ of
invariants is commutative. In this paper we classify all multiplicity
free symplectic representations.

}
\bigskip

\beginsection Introduction. Introduction

Let $V$ be a finite dimensional complex vector space equipped with a
symplectic form $\omega$. Then we can form the {\it Weyl algebra}
$\cW(V)$ which, by definition, is generated by $V$ with relations
$$
v\cdot w-w\cdot v=\omega(v,w).
$$
Observe that the symplectic group $Sp(V)$ acts on $\cW(V)$ by
algebra automorphisms.

Assume a connected reductive group $G$ is acting linearly on $V$
preserving the symplectic structure $\omega$. Then $V$ is called a
{\it symplectic representation} of $G$. One can think of a symplectic
representation as a homomorphism $\rho:G\pfeil Sp(V)$. Hence $G$ will
also act on the Weyl algebra $\cW(V)$.

\Definition: The symplectic representation $G\pfeil Sp(V)$ is called
{\it multiplicity free\/} (or an ``MFSR'') if $\cW(V)^G$ is commutative.
\medskip
\noindent Another way to phrase the definition is as follows. Recall
(see, e.g., \cite{Howe}) that there is a Lie algebra homomorphism
$\sp(V)\into\cW(V)$ whose adjoint action integrates to the
$Sp(V)$\_action. Thus we get a Lie algebra homomorphism
$$8
\tilde\rho:\fg\pfeil\cW(V)
$$
such that $\cW(V)^G$ is simply the commutant of $\tilde\rho(\fg)$ in
$\cW(V)$.

There are two instances of multiplicity free symplectic
representations which have been studied extensively. The first one are
Howe's reductive dual pairs \cite{Howe}. More precisely, let
$$9
G=G_1\times G_2=Sp_{2m}(\CC)\times SO_n(\CC)\hbox{ acting on }
V=\CC^{2m}\otimes\CC^n\quad{\rm or}
$$
$$7
G=G_1\times G_2=GL_m(\CC)\times GL_n(\CC)\hbox{ acting on }
V=(\CC^m\otimes\CC^n)\oplus(\CC^m\otimes\CC^n)^*
$$
Then Howe showed that the algebras generated by $\fg_1$ and $\fg_2$
are mutual commutants\footnote{In situation \cite{E9} this is not
quite true. But still the commutant of $\fg_1$ is the
algebra generated by $\fg_2$ which suffices for our purposes.} inside
$\cW(V)$. In particular, the commutant of $\fg$ is commutative.

The other class of multiplicity free symplectic representations
studied previously generalizes the case \cite{E7} above. More
precisely, let $G\pfeil GL(U)$ be any finite dimensional
representation. Then $V=U\oplus U^*$ carries the $G$\_invariant
symplectic structure
$$6
\omega(u_1+u_1^*,u_2+u_2^*)=\<u_1^*,u_2\>-\<u_2^*,u_1\>.
$$
In this situation one can identify $\cW(V)$ with the algebra
$\cP\cD(U)$ of polynomial coefficient linear differential operators on
$U$. From that fact one can deduce that $(G,V)$ is multiplicity free
if and only if the algebra $\CC[U]$ of polynomial functions on $U$ is
multiplicity free, i.e., does not contain any simple $G$\_module more
than once. Modules $U$ of this type are also called {\it multiplicity
free spaces}. They have been thoroughly investigated, e.g., in
\cite{HU}. In particular, they have been classified by Kac~\cite{Kac},
Benson\_Ratcliff~\cite{BR}, and Leahy~\cite{Leahy}. In the present
paper we extend this work by classifying all multiplicity free
symplectic representations

We continue by stating some fundamental properties of multiplicity
free symplectic representations, proved in \cite{inprep}. Let
$\cZ(\fg)$ be the center of the universal enveloping algebra
$\cU(\fg)$. Then the homomorphism \cite{E8} induces algebra
homomorphisms
$$
\cU(\tilde\rho):\cU(\fg)\pfeil\cW(V)\quad{\rm and}
\quad\cU(\tilde\rho)^G:\cZ(\fg)\pfeil\cW(V)^G.
$$
Clearly, if $\cU(\tilde\rho)^G$ is surjective then $V$ is multiplicity
free. In general, one can show that $V$ is a MFSR if and only if
$\cU(\tilde\rho)^G$ is ``almost'' surjective, i.e., $\cW(V)^G$ is a
finitely generated $\cZ(\fg)$\_module.

The Weyl algebra is filtered by putting $V$ in degree $1$. The
associated graded version of $\cU(\tilde\rho)$ is the {\it moment
map}\footnote{Actually, $m$ is the composition of
$\|gr|\cU(\tilde\rho)$ with the isomorphism $V\pfeil V^*$ induced by
$\omega$.}
$$
m:V\pfeil\fg^*:v\mapsto[\xi\mapsto\textstyle{1\over2}\omega(\xi v,v)].
$$
One can show that $V$ is a MFSR if and only if ``almost'' all
$G$\_invariants on $V$ are pull\_backs of coadjoint invariants, i.e.,
$\cO(V)^G$ is a finitely generated $\cO(\fg^*)^G$\_module. 

The central result about multiplicity free symplectic representation
is the following theorem proved in \cite{inprep}:

\Theorem cofree. Let $G\pfeil Sp(V)$ be multiplicity free symplectic
representation. Then $V$ is cofree, i.e., the algebra of
invariants $\cO(V)^G$ is a polynomial ring and $\cO(V)$ is a
free $\cO(V)^G$\_module.

\noindent Observe that this has an immediate corollary, namely that
also the algebra $\cW(V)^G$ is a polynomial ring and that $\cW(V)$ is a
free $\cW(V)^G$\_module (left or right).

The identification of $\cO(V)^G$ with a polynomial ring can be made
more precise. Recall that $\cO(\fg^*)^G$ can be identified with
$\cO(\ft^*)^W$ where $\ft\subseteq\fg$ is a Cartan subalgebra and $W$
is the Weyl group. Then one can find a subspace $\fa^*\subseteq\ft^*$
and a subgroup $W_V\subseteq N_W(\fa^*)/C_W(\fa^*)$ with
$\cO(V)^G\cong\cO(\fa^*)^{W_V}$ and such that
$$
\matrix{
V&\pf m&\fg^*\cr
\untenPf&&\untenPf\cr
\fa^*/W_V&\pfeil&\ft^*/W\cr}
$$
commutes. Moreover, $W_V$ is a reflection group.

\beginsection Classif. Statement of the classification 

Before we state the result of our classification we have to develop
some terminology. For ease of notation we state all of our results in
terms of the reductive Lie algebra $\fg$. The $n$\_dimensional
commutative Lie algebra will be denoted by $\ft^n$. We will only
consider algebraic representations, i.e., those which come from a
representation of the corresponding group $G$. This is only an issue
if $\fg$ is not semisimple.

First, we focus on symplectic representations of a fixed Lie algebra
$\fg$.

\Definition: a) A symplectic representation is called {\it
indecomposable} if it is not isomorphic to the sum of two non\_trivial
symplectic representations.

b) Let $V$ be a symplectic representation. Then $V$ is said to be of
{\it type~1} if $V$ is irreducible as a $\fg$\_module. It is of {\it
type~2} if $V=U\oplus U^*$ where $U$ is an irreducible $\fg$\_module
not admitting a symplectic structure and the form is \cite{E6}.

\medskip

\Theorem. \item{a)} Every indecomposable symplectic representation is either
of type~1 or~2.\Par
 \item{b)} Assume two symplectic representation are isomorphic as
$\fg$\_modules. Then they are isomorphic as symplectic representations.\Par
 \item{c)} Every symplectic representation is a direct sum of finitely many
indecomposable symplectic representations. The summands are unique up
to permutation.\Par

\Proof: Let $V$ be a symplectic representation and $U\subseteq V$ an
irreducible submodule. If $\omega|_U\ne0$ then $V=U\oplus U^\perp$. In
particular, $V$ has a type~1 summand. If $\omega|_U=0$ then $U$ is in
the kernel of $V\pfeil U^*:v\mapsto\omega(v,\cdot)$. By complete
reducibility, there is an inclusion $U^*\into V$ such that $\omega$
restricted to $\Uq:=U\oplus U^*$ is the standard form. Any
symplectic structure $\omega_0$ on $U$ induces an isomorphism $U\pfeil
U^*:u\mapsto u^*$. Moreover,
$$
\omega(u_1+u_1^*,u_2+u_2^*)=\<u_1^*,u_2\>-\<u_2^*,u_1\>=
\omega_0(u_1,u_2)-\omega_0(u_2,u_1)=2\omega_0(u_1,u_2).
$$
This means that $\omega$ restricted to the diagonal sitting in
$U\oplus U\cong\Uq$ is non\_zero bringing us back to the first
case. Therefore, we may assume that $U$ does not admit a symplectic
structure. Since $V=\Uq\oplus\Uq^\perp$ we see that $V$ has a type~2
summand. This already proves part a).

Let $(\tilde V,\tilde \omega)$ be a second symplectic representation
and assume $V\cong\tilde V$ as $\fg$\_module. Then $\tilde V$ has a
summand $\tilde U$ isomorphic to $U$. Assume first that $U$ is
symplectic. Then we can choose $\tilde U$ such that $\tilde
\omega|_{\tilde U}\ne0$. Since the symplectic structure of an
irreducible module is unique up to a scalar we can find an symplectic
isomorphism $U\pfeil\tilde U$. The orthogonal spaces $U^\perp\cong
V/U$ and $\tilde U^\perp=\tilde V/\tilde U$ are isomorphic as
$\fg$\_modules. By induction, they are isomorphic as symplectic
representations. Thus we get a symplectic isomorphism
$$3
V=U\oplus U^\perp\pf\sim\tilde U\oplus\tilde U^\perp=\tilde V.
$$
The same argument works if $U$ does not have a symplectic
structure; just replace $U$ by $\Uq$ in \cite{E3}. This completes the
proof of b). The preceding discussion also shows how to read off the
components of a symplectic representation from its decomposition as
$\fg$\_module. This proves c).\qed

\noindent Because of part b) we won't need to explicitly specify the
symplectic structure in our tables.

Now we set up notation for when $\fg$ varies.

\Definition: Let $\rho_1:\fg_1\pfeil\sp(V_1)$,
$\rho_2:\fg_2\pfeil\sp(V_2)$ be two symplectic representations.

\item{a)} $V_1$ and $V_2$ are {\it equivalent} if there is an
symplectic isomorphism $\phi:V_1\pfeil V_2$ (inducing an isomorphism
$\phi:\sp(V_1)\pfeil\sp(V_2)$) such that
$\rho_2(\fg_2)=\phi(\rho_1(\fg_1))$.

\item{b)} The {\it product} of $V_1$ and $V_2$ is the algebra $\fg_1+\fg_2$
acting on $V_1\oplus V_2$. A symplectic representation is called {\it
connected} if it is not equivalent to the product of two non\_trivial
symplectic representations.

\medskip

\noindent This definition of equivalence has two consequences which
have to be kept in mind: first, it depends only on the image of $\fg$,
i.e., the kernel is being ignored. Second, two representations which
differ by an (outer) automorphism are equivalent. This holds, for
example, for the two spin representations of $\so_{2n}$.

Two symplectic representations are multiplicity free if and only if
their product is. This follows, e.g., directly from the
definition. Therefore, it suffices to classify connected
representations. Unfortunately, a representation may be connected for
a rather trivial reason. Take, e.g., finitely many
(non\_symplectic) representations $(\fg_1,U_1),\ldots,(\fg_s,U_s)$ and
form the symplectic representation of $\fg=\prod_i\fg_i$ acting on
$V=\oplus_i(U_i\oplus U_i^*)$. The one\_dimensional Lie algebra
$\ft^1$ acts on each summand by $t\cdot(u,u^*)=(tu,-tu^*)$. This way,
we get a connected representation of $\fg+\ft^1$. There are cases
where $(\fg+\ft^1,V)$ is multiplicity free while $(\fg,V)$ is
not. Thus, we cannot simply ignore the center of $\fg$. Instead, we go
up and enlarge the algebra $\fg+\ft^1$ to $\fg+\ft^s$. Then the
representation becomes disconnected namely the product of the
$(\fg+\ft^1,U_i\oplus U_i^*)$.

\Definition: A symplectic representation $\rho:\fg\pfeil\sp(V)$ is
{\it saturated} if $\rho(\fg)$ is its own normalizer in $\sp(V)$.

\medskip

\noindent Every type~2 representation $U\oplus U^*$ has non\_trivial
endomorphisms namely $\ft^1$ acting by
$t\cdot(u,u^*)=(tu,-tu^*)$. Roughly speaking, saturatedness means that
$\rho(\fg)$ contains all these endomorphisms. More precisely:

\Proposition satcrit. A symplectic representation
$\rho:\fg\pfeil\sp(V)$ is saturated if and only if every type~1
component appears with multiplicity one and the number of type~2
components equals the dimension of the center of $\rho(\fg)$.

\Proof: Let $V=\sum_iC_i^{n_i}$ be the decomposition of $V$ into
components such that the $C_i$ are pairwise non\_isomorphic. The
centralizer of $\rho(\fg)$ in $\sp(V)$ is the product of the
centralizers of the $C_i^{n_i}$. There are three cases to consider:

1. $C_i$ is of type 1. Then the centralizer is $\so_{n_i}$. 

2a. $C_i=U\oplus U^*$ is of type 2 with $U\not\cong U^*$. Then the
   centralizer is $\gl_{n_i}$.

2b. $C_i=U\oplus U^*$ is of type 2 with $U\cong U^*$. Then the
   centralizer is $\sp_{2n_i}$.

\noindent
Now assume that $V$ is saturated. Then $\rho(\fg)$ contains its
centralizer. Therefore, $n_i=1$ in all cases, there are no components
of type~2b and the dimension of the center equals the number of
components of type 2.

Conversely, if $n_i=1$ for type~1 components then the dimension of the
center of $\rho(\fg)$ is at most the number of the $C_i$ of
type~2a. This implies that there are no type~2b components and all
$n_i=1$. Moreover, $\rho(\fg)$ will contain its centralizer, i.e., is
self\_normalizing.\qed

\noindent This proposition shows in particular that, up to
equivalence, a saturated representation can be easily reconstructed
from the representation of the semisimple part $\fg'$ of $\fg$. More
precisely, the type~1 components of $(\fg,V)$ are those symplectic
irreducible summands of $(\fg',V)$ which occur with odd
multiplicity. The rest of the irreducible summands can be paired up in
the form $M\oplus M^*$. These pairs form the type~2 components of
$(\fg,V)$. Moreover, $\fg=\fg'\oplus\ft^s$ where $s$ is the number of
these pairs. For example, $(\fg',V)=(\sp_n,(\CC^n)^{\oplus7})$
corresponds to the saturated representation of $\fg=\sp_n\oplus\ft^3$
on $V=\CC^n\oplus(\CC^n\oplus\CC^n)^{\oplus3}$.

For that reason, we use in the sequel the following notation: let
$\rho:\fs\pfeil\gl(U)$ be a representation of a semisimple algebra
$\fs$. Then we denote the type~2 representation of $\fg=\fs+\ft^1$ on
$U\oplus U^*$ by $T(U)$. Continuing, if $U_1,U_2$ are two
representations of $\fs$ then $T(U_1)\oplus T(U_2)$ is a
representation of $\fg=\fs+\ft^2$. Finally, one should always keep in
mind that $T(U)$ is equivalent to $T(U^*)$.

The classification of multiplicity free symplectic
representations is easily reduced to the saturated case:

\Theorem saturiert. Let $\rho:\fg\pfeil\sp(V)$ be a multiplicity free
symplectic representation. Then there is a unique reductive subalgebra
$\overline\fg\subseteq\sp(V)$ with
$\overline\fg'\subseteq\rho(\fg)\subseteq\overline\fg$ such that
$(\overline\fg,V)$ is saturated and multiplicity free.\Par Conversely,
let $(\overline\fg, V)$ be a saturated and multiplicity free symplectic
representation. Let $\fc$ be the center of $\overline\fg$. Then there
is a unique subspace $\fa\subseteq\fc$ with the following
property:\Par Let $\fd\subseteq\fc$ be a subspace and
$\fg:=\overline\fg'\oplus\fd$. Then $(\fg,V)$ is multiplicity free if
and only if $\fd+\fa=\fc$.

\noindent The proof of this and the next three theorems will be given
in section~\cite{Beweis}.

Now, we are in the position to state our classification results. First
the indecomposable ones:

\Theorem IndecompClass. All saturated indecomposable MFSRs are listed
in Table~1 (type~1) and Table~2 (type~2).

\noindent Here, and throughout the paper we are using the convenient
notation that a Lie algebra denotes also its defining
representation. In case of an ambiguity (e.g., $\spin_{2n}$ and
$\sE6$) one may choose either of the two candidates. For example,
$\sp_{2m}\otimes\spin_7$ is the representation of $\fg=\sp_{2m}+\so_7$
on $V=\CC^{2m}\otimes\CC^8$ where $\so_7$ acts on $\CC^8$ via its
spin\_representation. The third fundamental representation of $\sp_6$
is denoted by $\wedge^3_0\sp_6$.

\midinsert
\def\tablename{1}\setbox0\hbox{$\langle\hbox{\tablename.}00\rangle$ }\SpBr=\wd0 
\centerline{{\bf Table 1:} MFSRs of type 1}
\smallskip
\centerline{\boxit{
\halign{\weiter$#$\hfill\quad&$#$\hfill\quad&$#$\hfill\quad&\hideWeyl{$#$\quad}\hfill
&$#$\hfill\quad&\hideIndex{\quad$#$}\hfill\cr
\nweiter(\fg,V)&&\hbox{rank}&W_V&\fl&i\cr
\noalign{\smallskip\hrule\smallskip}
\fs\fp_{2m}\otimes\fs\fo_p&2m\ge p=2n\ge4&n&\sD n&\fs\fp_{2m-p}+\ft^n
&1+\cr
\nweiter&2m\ge p=2n+1\ge3&n&\sB n
&\fs\fp_{2m-p}+\ft^n&1+\cr
\nweiter&2\le2m<p&m&\sC m
&\fs\fo_{p-2m}+\ft^m&1+\cr
\fs\fp_{2m}&m\ge1&0&\Nichts
&\fs\fp_{2m-1}&1+\cr
\noalign{\smallskip\hrule\smallskip}
\fs\fp_{2m}\otimes\spin_7&m=1&1&\sA1&\fs\fl_3+\ft^1&
1+\cr
\nweiter&m=2&3&\sC2+\sA1&\ft^2&
1+\cr
\nweiter&m=3&6&\sC3+\sB3&0&
1+\cr
\nweiter&m\ge4&7&\sD4+\sB3
&\fs\fp_{2m-8}&2+\cr
\fs\fl_2\otimes\spin_9&&2&2\sA1&\fs\fl_3+\ft^1&
1+\cr
\spin_n&n=11&1&\sA1&\fs\fl_5&
1+\cr
\nweiter&n=12&1&\sA1&\fs\fl_6&1+\cr
\nweiter&n=13&2&\sB2&\fs\fl_3+\fs\fl_3&
1+\cr
\noalign{\smallskip\hrule\smallskip}
S^3\fs\fl_2&&1&\sA1&0&
1+\cr
\wedge^3\fs\fl_6&&1&\sA1&\fs\fl_3+\fs\fl_3&
1+\cr
\wedge^3_0\fs\fp_6&&1&\sA1&\fs\fl_2&
1+\cr
\fs\fl_2\otimes\sG2&&1&\sA1&\fs\fl_2+\ft^1&
1+\cr
\fs\fp_4\otimes\sG2&&4&\sC2+\sG2&0&
1+\cr
\sE7&&1&\sA1&\sE6&1+\cr
}}}
\endinsert

\noindent The tables list also the rank of $V$ (i.e., the dimension of
$\CC[V]^\fg$) and the generic isotropy algebra $\fl$. To reduce the
number of cases we define $\gl_0=\so_1=\sp_0=\sp_{-1}:=0$. Moreover,
$\sp_{2n-1}$ denotes the isotropy algebra of any non\_zero vector in
the defining representation of $\sp_{2n}$. In the last column ``$i$''
we list some properties of the homomorphism
$\overline\phi^G:\cO(\fg^*)^G\pfeil\cO(V)^G$, the graded version of
$\phi^G:\cZ(\fg)\pfeil\cW(V)^G$, which might be of interest. This
homomorphism factors as
$$
\cO(\fg^*)^G=\cO(\ft^*)^W\pf\alpha\cO(\fa^*)^N
{\buildrel\beta\over\into}\cO(V)^G
$$
where $\fa^*\subseteq\ft^*$ is a certain subspace (namely the span of
$\Phi_+^t$ from the algorithm of section~\cite{Tools} below) and $N$
is its normalizer in the Weyl group $W$. The group $N$ acts on $\fa^*$
always as a reflection group even though there is no a priori
reason. Therefore, $\cO(\fa^*)^N$ is a polynomial ring and $\cO(V)^G$
is a free $\cO(\fa^*)^N$\_module. The last column records its rank
$[N:W_V]$. The map $\alpha$ is almost surjective in the sense that
$\|Image|\alpha$ and $\cO(\fa^*)^N$ have the same field of
fractions. The sign behind the rank signifies whether $\alpha$ is
surjective or not. This means in particular that $\overline\phi^G$ is
surjective if and only if the last column contains a $1+$. It is
injective, if and only if $\fl=0$.

In the classification of decomposable representation it turns out that
certain $\sl_2$\_factors in $\fg$ pose complications.

\midinsert
\xxx0\def\tablename{2}\setbox0\hbox{$\langle\hbox{\tablename.}00\rangle$ }\SpBr=\wd0 
\centerline{{\bf Table 2:} MFSRs of type 2}
\smallskip
\centerline{%
\boxit{\halign{\weiter$T(#)$\hfill\quad&$#$\hfill\quad&$#$\hfill\quad&\hideWeyl{$#$\quad}\hfill&$#$\hfill\quad&\hideIndex{\quad$#$}\cr
\omit\hfill$(\fg,V)$\hfill&&\hbox{rank}\hskip-5pt&W_V&\fl&i\cr
\noalign{\smallskip\hrule\smallskip}
\fs\fl_m\otimes\fs\fl_n&m\ge n\ge2&n&\sA{n-1}&\fg\fl_{m-n}+\ft^{n-1}&1+\cr
\wedge^2\fs\fl_n&n=2m\ge4&m&\sA{m-1}&\fs\fl_2^m&
1+\cr
\nweiter&n=2m+1\ge5&m&\sA{m-1}&\fs\fl_2^m+\ft^1&
1+\cr
S^2\fs\fl_n&n\ge2&n&\sA{n-1}&0&
1+\cr
\fs\fl_n&n\ge2&1&\Nichts&\fg\fl_{n-1}&1+\cr
\noalign{\smallskip\hrule\smallskip}
\fs\fp_{2m}&m\ge2&1&\Nichts&\fs\fp_{2m-2}+\ft^1&
1+\cr
\fs\fp_{2m}\otimes\fs\fl_n&m\ge2,n=2&3&2\sA1&\fs\fp_{2m-4}+\ft^1&
1+\cr
\nweiter&m=2,n=3&5&\sC2+\sA2&0&
1+\cr
\nweiter&m=2,n\ge4&6&\sC2+\sA3&\fg\fl_{n-4}&
1+\cr
\nweiter&m\ge3,n=3&6&\sA3+\sA2&\fs\fp_{2m-6}&
2+\cr
\noalign{\smallskip\hrule\smallskip}
\fs\fo_m&m\ge5&2&\sA1&\fs\fo_{m-2}&1+\cr
\spin_n&n=7&2&\sA1&\fs\fl_3&
1+\cr
\omit\hfill&n=9&3&2\sA1&\fs\fl_3&
1-\cr
\omit\hfill&n=10&2&\sA1&\fs\fl_4+\ft^1&
1+\cr
\noalign{\smallskip\hrule\smallskip}
\sG2&&2&\sA1&\fs\fl_2&
1+\cr
\sE6&&3&\sA2&\fs\fo_8&1-\cr
}}}
\endinsert

\Definition: A MFSR $(\fg,V)$ is said to have an {\it
$\sl_2$\_link} if $\fg$ has a factor which is isomorphic to $\sl_2$
and which acts effectively on more than one component of $V$.

\Theorem CompClass. All connected saturated MFSRs without
$\sl_2$\_links are listed in tables~11,~12, and~22.

\noindent In these tables we used a notation which is best explained
by way of an example:
$$
\sl_2\otimes\so_7\vplus\spin_7\otimes\sl_2
$$
means that $\fg=\sl_2+\so_7+\sl_2$ is acting on
$V=\CC^2\otimes\CC^7\oplus\CC^8\otimes\CC^2$ where the first/last
$\sl_2$ is acting on the first/last $\CC^2$ only while $\so_7$ is acting
diagonally on $\CC^7$ (defining representation) and $\CC^8$ (spin
representation). The line under the $\oplus$\_sign means that the
algebras immediately to the left and to the right are being identified
and is acting diagonally.

\midinsert
\xxx0\def\tablename{11}\setbox0\hbox{$\langle\hbox{\tablename.}00\rangle$ }\SpBr=\wd0 
\centerline{{\bf Table 11:} MFSRs of type 1--1 without $\sl_2$\_links}
\smallskip
\centerline{%
\boxit{{\def\vplus{&}
\halign{\weiter\hfill$#\ \underline\oplus\ $&$#$\hfill&$#$\hfill\ &$#$\hfill\quad&\hideWeyl{$#$\quad}\hfill&$#$\hfill\hideEmb{\quad}&\hideIndex{\quad$#$}\cr
\omit\hfill$(\fg,V)$&&&\hbox{rank}&W_V&\fl&i\cr
\noalign{\smallskip\hrule\smallskip}
\sl_2\otimes\so_n\vplus\so_n\otimes\sl_2
&n\ge7&4&2\sA1+\sB2&\so_{n-4}\subset\so_n&1+\cr
\spin_{12}^+\vplus\spin_{12}^-
&&4&2\sA1+\sB2&\sl_2+\sl_2&1-\cr
\sl_2\otimes\so_{12}\vplus\spin_{12}
&&3&3\sA1&\sl_4+\ft^1&
1-\cr
\sp_4\otimes\so_{12}\vplus\spin_{12}
&&7&\sC2{+}\sA1{+}\sD4&\sl_2&1-\cr
\sl_2\otimes\so_{11}\vplus\spin_{11}
&&4&2\sA1+\sB2&\sl_3\subset\so_{11}&1-\cr
\sl_2\otimes\so_8\vplus\spin_8\otimes\sl_2
&&3&3\sA1&\sl_2+\ft^2&
1+\cr
\sp_4\otimes\so_8\vplus\spin_8\otimes\sl_2
&&7&\sC2{+}\sD4{+}\sA1&0&1+\cr
\sl_2\otimes\so_7\vplus \spin_7\otimes\sl_2
&&4&2\sA1+\sB2&\ft^1\subset0{+}\so_7{+}\sl_2&
1+\cr
\sl_2\otimes\spin_7\vplus \spin_7\otimes\sl_2
&&5&3\sA1+\sB2&0&3+\cr
\noalign{\smallskip\hrule\smallskip}
\sl_2\otimes\so_6\vplus\so_6\otimes\sl_2
&&4&2\sA1+\sB2&\ft^1\subset\so_6&1+\cr
\noalign{\smallskip\hrule\smallskip}
\so_p\otimes\sp_{2m}\vplus\sp_{2m}
&\hideWeyl{\kern-20pt} p>2m\ge2&2m&\sB m+\sC m&
\so_{p-2m}&1+\cr
\omit\hfill&\omit\hfill&\hideWeyl{\kern-45pt} 3\le p\hideWeyl{=2n}\le2m&p&\sD n+\sC n&\sp_{2m-p-1}&1+\cr
\hideWeyl{\omit\hfill&\omit\hfill&
\kern-50pt3\le p=2n{-}1{<}2m
&p&\sB{n-1}+\sD n&\sp_{2m-p-1}&1+\cr}
\wedge^3_0\sp_6\vplus\sp_6
&&2&2\sA1&\sl_2&1-\cr
\spin_7\otimes\sp_4\vplus\sp_4
&&5&\sA1{+}\sB2{+}\sC2&0&3+\cr
\sp_{2m}\otimes\so_5\vplus\sp_4
&m=1&2&2\sA1&\ft^1&
1+\cr
\omit\hfill&\omit\hfill
&m\ge2&4&2\sC2&\sp_{2m-5}&
1+\cr
\sl_2\otimes\so_5\vplus\so_5\otimes\sl_2
&&4&2\sA1+\sB2&0&1+\cr
}}}}
\endinsert
\midinsert
\xxx0\def\tablename{12}\setbox0\hbox{$\langle\hbox{\tablename.}00\rangle$ }\SpBr=\wd0 
\centerline{{\bf Table 12:} MFSRs of type 1--2 without $\sl_2$\_links}
\smallskip
\centerline{%
\boxit{
{\def\vplus{&}
\halign{\weiter\hfill$#\ \underline\oplus\ $&$#$\hfill&$#$\hfill\ &$#$\hfill\quad&\hideWeyl{$#$\quad}\hfill&$#$\hfill\hideEmb{\quad}&\hideIndex{\quad$#$}\cr
\omit\hfill$(\fg,V)$&&&\hbox{rank}&W_V&\fl&i\cr
\noalign{\smallskip\hrule\smallskip}
\spin_{12}\vplus T(\so_{12})
&&4&3\sA1&\sl_4&2-\cr
\sl_2\otimes\so_{10}\vplus T(\spin_{10})
&&6&2\sA1+\sA3&\sl_2&1-\cr
\sl_2\otimes\so_8\vplus T(\spin_8)
&&4&3\sA1&\sl_2+\ft^1&
1-\cr
\sl_2\otimes\so_7\vplus T(\spin_7)
&&5&2\sA1+\sB2&0&3+\cr
\noalign{\smallskip\hrule\smallskip}
\wedge^3\sl_6\vplus T(\sl_6\otimes\sl_2)
&&6&2\sA1+\sA3&\ft^1&
1-\cr
\wedge^3\sl_6\vplus T(\sl_6)
&&3&2\sA1&\sl_2+\sl_2+\ft^1&1-\cr
\sp_{2m}\otimes\so_6\vplus T(\sl_4)
&m=1&3&2\sA1&\ft^2&
1+\cr
\omit\hfill&\omit\hfill
&m=2&6&\sC2+\sA3&0&1+\cr
\omit\hfill&\omit\hfill
&m\ge3&7&2\sA3&\sp_{2m-6}&2+\cr
\sl_2\otimes\so_6\vplus T(\sl_4\otimes\sl_2)
&&6&2\sA1+\sA3&0&1+\cr
\noalign{\smallskip\hrule\smallskip}
\sp_{2m}\vplus T(\sp_{2m})
&m\ge2&2&\sA1&\sp_{2m-3}&1+\cr
\wedge^3_0\sp_6\vplus T(\sp_6)
&&4&\sA1+\sC2&0&3+\cr
\sp_4\vplus T(\so_5)&&3&2\sA1&0&2+\cr
\sl_2\otimes\so_5\vplus T(\sp_4)
&&4&\sA1+\sB2&0&1+\cr
}}}}
\endinsert
\midinsert
\xxx0\def\tablename{22}\setbox0\hbox{$\langle\hbox{\tablename.}00\rangle$ }\SpBr=\wd0 
\centerline{{\bf Table 22:} MFSRs of type 2--2 without $\sl_2$\_links}
\smallskip
\centerline{%
\boxit{
{\def\vplus{&}
\halign{\weiter\hfill$#\ \underline\oplus\ $&$#$\hfill&$#$\hfill\ &$#$\hfill\quad&\hideWeyl{$#$\quad}\hfill&$#$\hfill\hideEmb{\quad}&\hideIndex{\quad$#$}\cr
\omit\hfill$(\fg,V)$&&&\hbox{rank}&W_V&\fl&i\cr
\noalign{\smallskip\hrule\smallskip}
T(\so_8)\vplus T(\spin_8)
&&5&3\sA1&\sl_2&1-\cr
\noalign{\smallskip\hrule\smallskip}
T(\wedge^2\sl_n)\vplus T(\sl_n)&
\hideWeyl{\kern-5pt} n\hideWeyl{=2m}\ge4&n&2\sA{m-1}&\ft^1&1-\cr
\hideWeyl{\omit\hfill&\omit\hfill&
\kern-22pt n=2m+1\ge 5&n&\sA{m}+\sA{m-1}&\ft^1&1-\cr}
T(\sl_m\otimes\sl_n)\vplus T(\sl_n)
&m\ge n\ge3&2n&2\sA{n-1}&\gl_{m-n}&1+\cr
\omit\hfill&\omit\hfill
&2\le m<n&2m+1&\sA{m}+\sA{m-1}&\gl_{n-m-1}&1+\cr
T(\sl_n)\vplus T(\sl_n)
&n\ge3&3&\sA1&\gl_{n-2}&1-\cr
\noalign{\smallskip\hrule\smallskip}
T(\sp_{2m})\vplus T(\sp_{2m})
&m\ge2&4&2\sA1&\sp_{2m-4}&2+\cr
}}}}
\endinsert

Representations with $\sl_2$\_links are dealt with in the next theorem.

\Theorem sl2link. All connected saturated MFSRs with $\sl_2$\_links
are obtained by taking any collection of representations from Table~S
and identifying any number of disjoint pairs of underlined
$\sl_2$'s. Moreover, not allowed is the identification of the two
$\sl_2$'s of {\rm(S.1)} and the combination {\rm(S.9)}+{\rm(S.9)}.

\midinsert
\xxx0\def\tablename{S}\setbox0\hbox{$\langle\hbox{\tablename.}00\rangle$ }\SpBr=\wd0 
\centerline{{\bf Table S}}
\smallskip
\centerline{%
\boxit{\halign{\weiter$#$\hfill\quad&$#$\hfill\quad&$#$\hfill\quad&$#$\hfill\hideEmb{\quad}\cr
\nweiter\fg,X&&\hbox{rank}&\fl\cr
\noalign{\smallskip\hrule\smallskip}
\usl_2\otimes\sp_{2m}\otimes\usl_2&m=1&1&\ft^2\cr
\nweiter&m\ge 2&2
&\sp_{2m-4}+\ft^2
\cr
\usl_2\otimes\so_8\vplus\spin_8\otimes\usl_2
&&3&\sl_2+\ft^2\cr
\noalign{\smallskip\hrule height .8pt\smallskip}
\so_n\otimes\usl_2&n\ge3&1
&\so_{n-2}+\ft^1\cr
\spin_{12}\vplus\so_{12}\otimes\usl_2
&&3&\sl_4+\ft^1\cr
\spin_9\otimes\usl_2&&2&\sl_3+\ft^1\cr
T(\so_8)\vplus\spin_8\otimes\usl_2
&&4&\sl_2+\ft^1\cr
\spin_7\otimes\usl_2&&1&\sl_3+\ft^1\cr
\sl_2\otimes\so_7\vplus\spin_7\otimes\usl_2
&&4&\ft^1\cr
\noalign{\smallskip\hrule\smallskip}
\usl_2&&0&\sp_1\cr
T(\usl_2)&&1&\ft^1\cr
T(\sl_m\otimes\usl_2)&m\ge2&2&\gl_{m-2}+\ft^1
\cr
T(\sl_4)\vplus\so_6\otimes\usl_2
&&3&\ft^2\cr
\noalign{\smallskip\hrule\smallskip}
\sp_{2m}\otimes  S^2\usl_2&m\ge1&1
&\sp_{2m-3}+\ft^1\cr
T(\sp_{2m}\otimes\usl_2)&m\ge2&3&\sp_{2m-4}+\ft^1\cr
\sp_4\vplus\so_5\otimes\usl_2
&&2&\ft^1\cr
\noalign{\smallskip\hrule\smallskip}
\sG2\otimes\usl_2&&1&\sl_2+\ft^1\cr
}}}
\endinsert

\noindent For example, if we identify the underlined $\sl_2$'s of
(S.6) and (S.8) then we get the algebra
$$
\fg=\ft^1+\so_8+\sl_2+\so_7+\sl_2
$$
acting on 
$$
V=(\CC^8\oplus\CC^8)\oplus(\CC^8\otimes\CC^2)\oplus
(\CC^2\otimes\CC^8)\oplus(\CC^7\otimes\CC^2)
$$
It is easier to visualize this as a graph:
$$4
\vcenter{\hbox{\bdiadraw
\KlabelU(0,0)[0]{T}
\DstrichRO(0,0)
\labelO(1,1)[0]{\so_8}
\DstrichRU(1,1)
\PlabelL(2,0)[0]{\omega_4}
\DstrichRO(2,0)
\labelO(3,1)[0]{\sl_2}
\DstrichRU(3,1)
\PlabelR(4,0)[0]{\omega_3}
\DstrichRO(4,0)
\labelO(5,1)[0]{\so_7}
\DstrichRU(5,1)
\punkt(6,0)
\DstrichRO(6,0)
\labelO(7,1)[0]{\sl_2}
\ediadraw}}
$$
The entries (S.1) and (S.2) are special in that they possess two
underlined $\sl_2$'s. They can be used to build connected saturated
MFSRs with arbitrary many components. For example, if one inserts (S.1)
into \cite{E4} then one gets
$$
\vcenter{\hbox{\bdiadraw
\KlabelU(0,0)[0]{T}
\DstrichRO(0,0)
\labelO(1,1)[0]{\so_8}
\DstrichRU(1,1)
\PlabelL(2,0)[0]{\omega_4}
\DstrichRO(2,0)
\labelO(3,1)[0]{\sl_2}
\DstrichRU(3,1)
\punkt(4,0)
\labelO(4,1)[0]{\sp_{2m}}
\strichU(4,1)
\DstrichRO(4,0)
\labelO(5,1)[0]{\sl_2}
\DstrichRU(5,1)
\PlabelR(6,0)[0]{\omega_3}
\DstrichRO(6,0)
\labelO(7,1)[0]{\so_7}
\DstrichRU(7,1)
\punkt(8,0)
\DstrichRO(8,0)
\labelO(9,1)[0]{\sl_2}
\ediadraw}}
$$
Even a circular pattern is allowed:
$$
\vcenter{\hbox{%
\bdiadraw
\labelO(0,1)[0]{\sp_{2m}}
\strichU(0,1)
\punkt(0,0)
\DstrichRO(0,0)
\move(0 0)\rlvec(5 1)
\labelO(1,1)[0]{\sl_2}
\DstrichRU(1,1)
\punkt(2,0)
\DstrichRO(2,0)
\labelO(3,1)[0]{\so_8}
\DstrichRU(3,1)
\PlabelL(4,0)[0]{\omega_4}
\DstrichRO(4,0)
\labelO(5,1)[0]{\sl_2}
\ediadraw
}}
$$
Finally, note the two boundary cases
$$
\vcenter{\hbox{\bdiadraw
\PlabelL(0,0)[0]{\omega_4}
\labelO(0,1)[0]{\so_8}
\strichU(0,1)
\DstrichRU(0,1)
\DstrichRO(0,0)
\labelO(1,1)[0]{\sl_2}
\strichU(1,1)
\punkt(1,0)
\ediadraw}}
$$
and
$$
\vcenter{\hbox{\bdiadraw
\labelO(0,1)[0]{\sp_{2m}}
\DstrichRU(0,1)
\labelO(1,1)[0]{\sl_2}
\strichU(1,1)
\DstrichRU(1,1)
\punkt(1,0)
\DstrichRO(1,0)
\labelO(2,1)[0]{\sl_2}
\strichU(2,1)
\punkt(2,0)
\DstrichRO(2,0)
\labelO(3,1)[0]{\sp_{2n}}
\ediadraw}}\quad=\quad\sp_{2m}\otimes\so_4\vplus\so_4\otimes\sp_{2n}
$$
The representation $\sp_{2m}\otimes\so_3\vplus\so_3\otimes\sp_{2n}$
is multiplicity free, as well, since it is the concatenation of (S.13)
with itself.

\beginsection Tools. The tool box

Before we enter the proofs of the classification theorems we collect
a couple of tools. The most important one is an algorithm
which allows to decide whether a given symplectic representation is
multiplicity free. Moreover, it determines the rank and the generic
isotropy group.

The input of the algorithm are two subsets of the weight lattice $X$,
namely $\Delta$, the set of roots, and $\Phi$, the set of weights. The latter
is actually a multiset, i.e., each weight is counted with its
multiplicity.

\Definition:
\item{a)}$\chi\in\Phi$ is {\it extremal} if $\alpha\in\Delta$,
$\<\chi|\alpha^\vee\>>0$ implies $\chi+\alpha\not\in\Phi$.

\item{b)}$\chi\in\Phi$ is {\it toroidal} if
$\<\chi|\alpha^\vee\>=0$ for all $\alpha\in\Delta$.

\item{c)}An extremal weight $\chi$ is {\it singular} if
$2\chi\in\Delta$ and the multiplicity of $\chi$ is one.

\medskip

\noindent{\bf The algorithm:} {\it Let $\chi\in\Phi$ be an extremal weight
which is neither toroidal nor singular. If no such $\chi$ exists then
stop with output $(\Delta,\Phi)$. Otherwise, let
$P:=\{\alpha\in\Delta\mid \<\chi|\alpha^\vee\>>0\}$, $Q:=\chi-P$ and
perform the following replacements
$$5
\eqalign{
\Delta&\ersetze\Delta\setminus(P\cup-P)\cr
\Phi&\ersetze\Phi\setminus(Q\cup-Q).\cr}
$$}

Note that if $\chi$ were toroidal then $P=\emptyset$ and nothing would
happen. The extremality of $\chi$ ensures that
$\Delta\setminus(P\cup-P)$ is the root system of a Levi subalgebra
$\fl$. The non\_singularity of $\chi$ implies $Q\cap(-Q)=\emptyset$
(in the multiset sense). This ensures that at each step, $\Phi$ is the
set of weights of a symplectic $\fl$\_representation.

The algorithm finishes with a pair $(\Delta_0,\Phi_0)$. Then $\Phi_0$
is the disjoint union of the the toroidal weights $\Phi_0^t$ and the
singular weights $\Phi_0^s$. Since $-\Phi_0^t=\Phi_0^t$ (and $0$ has
even multiplicity) we can find a subset $\Phi_+^t$ such that
$\Phi_0^t$ is the disjoint union of $\Phi_+^t$ and $-\Phi_+^t$. Now,
our criterion is:

\Theorem. The symplectic representation is multiplicity free if and only if
the $\Phi_+^t\subseteq X$ is linearly independent. In that case, we
have
$$
\|dim|V\mod G=\|rk|V=|\Phi_+^t|=\textstyle{1\over2}|\Phi_0^t|.
$$

\noindent Also the generic isotropy group can be determined. Let
$L\subseteq G$ be the Levi subgroup with root system $\Delta_0$. Let
$L_0\subseteq L$ be the intersection of all $\|ker|\chi$ with
$\chi\in\Phi_0^t$. Call two singular weights $\chi_1,\chi_2$ equivalent
if $\chi_1-\chi_2\in\Delta_0$. Each equivalence class determines a
direct factor $\cong Sp_{2m}$ of $L_0$ where $2m$ is the size of the
class.

\Theorem. Let $L_0=L_1\times Sp_{2m_1}\times\ldots\times Sp_{2m_r}$ be the
decomposition determined by singular weight classes. Then the generic
isotropy group of $G$ in $V$ is (conjugate to) $L_1\times
Sp_{2m_1-1}\times\ldots\times Sp_{2m_r-1}$ (where $Sp_{2m-1}$ is an
isotropy group of $Sp_{2m}$ in $\CC^{2m}\setminus\{0\}$.

\medskip\noindent{\it Proofs:} Let $\chi\in\Phi$ be an non\_singular
extremal weight. Then $\chi$ is dominant with respect to a system
$\Delta^+\subseteq\Delta$ of positive roots. Let
$\Delta_+,\Delta_0,\Delta_-\subseteq\Delta$ be the set of roots
$\alpha$ with $\<\chi|\alpha^\vee\>>0,=0,<0$, and let
$\fp_u,\fl,\fp_u^-\subseteq \fg$ be the corresponding subalgebras. Let
$v_0\in V$ be a weight vector for $\chi$ and choose a weight vector
$v_0^-\in V$ whose weight is $-\chi$ with $\omega(v_0,v_0^-)=1$. Then
$V_0:=(\fp_u^- v_0)^\perp\cap(\fp_u v_0^-)^\perp$ is a symplectic
$\fl$\_representation. Let $\Phi_0$ be the set of weights of
$V_0$. Moreover, we have $\Delta_0=\Delta\setminus(P\cup-P)$ and
$\Phi_0=\Phi\setminus(Q\cup-Q)$. Thus, the next step in the algorithm
encodes the representation $(\fl,V_0)$.

It has been shown in \cite{inprep} that $(\fg,V)$ is multiplicity free
if and only if $(\fl,V_0)$ is. Moreover, the generic isotropy groups
are the same. Therefore, both properties can be read off of the output
$(\Delta_0,\Phi_0)$ of the algorithm which corresponds to a symplectic
representation of a Levi subgroup $\fl$ on a space $V_0$ on which
every weight is either toroidal or singular. If $\chi$ is a singular
dominant weight then $2\chi$ is a root. This means that $\chi$ belongs
to the defining representation of a subalgebra of type
$\sC m$. Therefore, we can decompose
$\fl=\fl_1+\sp_{2m_1}+\ldots+\sp_{2m_r}$ and
$V_0=V_1\oplus\CC^{2m_1}\oplus\ldots\oplus\CC^{2m_r}$ where $V_1$ is a
sum of one\_dimension $\fl_1$\_modules. This implies $V_1=U\oplus U^*$
where $\fl_1$ acts on $U$ with the characters
$\Phi_+^t$. Thus $(\fl,V_0)$ is multiplicity free if and only if
$(\fl_1,U)$ is a multiplicity free action if and only if $\Phi_+^t$ is
linearly independent. The claim about the generic isotropy group is
also obvious from this description.\qed

\Examples: 1. Let $G=\sl_6+\ft^2$ and $V=\wedge^3U\oplus(U\oplus
U^*)\oplus(U\oplus U^*)$ where $U=\CC^6$ is the defining
representation of $\sl_6$ and $\ft^2$ acts according to the pattern
$(0,s,-s,t,-t)$. Then
$$
\eqalign{
\Delta=&\,\{\epsilon_i-\epsilon_j\}\cr
\Phi=&\,\{\epsilon_i+\epsilon_j+\epsilon_k\}
\cup\{\pm(\epsilon_i+\eta)\}
\cup\{\pm(\epsilon_i+\eta')\}\cr}
$$
where $i,j,k=1,\ldots,6$ are pairwise distinct. Moreover, the relation
$\sum_i\epsilon_i=0$ holds. We start with the extremal
weight $\chi=\epsilon_1+\epsilon_2+\epsilon_3$. Then
$$
P=\{\epsilon_i-\epsilon_j\mid i=1,2,3;j=4,5,6\}
$$
and we are left with
$$
\eqalign{
\Delta=&\,\{\epsilon_i-\epsilon_j\mid 1\le i\ne j\le
3\}\cup\{\epsilon_i-\epsilon_j\mid 4\le i\ne j\le 6\}\cr
\Phi=&\,\{\pm(\epsilon_1+\epsilon_2+\epsilon_3)\}
\cup\{\pm(\epsilon_i+\eta)\}
\cup\{\pm(\epsilon_i+\eta')\}\cr}
$$
Now put successively $\chi=\epsilon_1+\eta$, $\chi=\epsilon_4+\eta$,
$\chi=\epsilon_2+\eta'$, and $\chi=\epsilon_5+\eta'$. The final result is
$$
\eqalign{
\Delta=&\,\emptyset\cr
\Phi=&\,\Phi_+^t\cup-\Phi_+^t,\Phi_+^t=\{
\epsilon_1+\epsilon_2+\epsilon_3,
\epsilon_1+\eta,
\epsilon_4+\eta,
\epsilon_1+\eta',
\epsilon_2+\eta',
\epsilon_4+\eta',
\epsilon_5+\eta'\}\cr}
$$
Since $\Phi_+^t$ is not linearly independent, the representation is
not multiplicity free.

2. The other example is $G=\sp_{2m}+\sp_4$ acting on
$V=\CC^{2m}\otimes\CC^5\oplus\CC^4$ where $\CC^5=\wedge^2_0\CC^4$ is
the second fundamental representation of $\sp_4(\CC)$. We assume
$m\ge2$. Then
$$
\eqalign{
\Delta=&\,\{\pm\epsilon_i\pm\epsilon_j,\pm2\epsilon_i\}
\cup\{\pm\epsilon_1'\pm\epsilon_2',\pm2\epsilon_1',\pm2\epsilon_2'\}\cr
\Phi=&\,\{\pm\epsilon_i\pm\epsilon_1'\pm\epsilon_2',\pm\epsilon_i\}
\cup\{\pm\epsilon_1',\pm\epsilon_2'\}\cr}
$$
We start with $\chi=\epsilon_1+\epsilon_1'+\epsilon_2'$. Then
$$
P=\{\epsilon_1\pm\epsilon_j\}\cup\{2\epsilon_1\}
\cup\{\epsilon_1'+\epsilon_2',2\epsilon_1',2\epsilon_2'\}
$$
and we get
$$
\eqalign{
\Delta=&\,\{\pm\epsilon_i\pm\epsilon_j,\pm2\epsilon_i\mid i,j\ge2\}
\cup\{\pm(\epsilon_1'-\epsilon_2')\}\cr
\Phi=&\,\{\pm(\epsilon_1+\epsilon_1'+\epsilon_2')\}
\cup\{\pm\epsilon_i\pm(\epsilon_1'-\epsilon_2'),\pm\epsilon_i\mid i\ge2\}
\cup\{\pm\epsilon_1',\pm\epsilon_2'\}\cr}
$$
Since $m\ge2$ we can repeat this process with
$\chi=\epsilon_2+\epsilon_1'-\epsilon_2'$ and get
$$
\eqalign{
\Delta=&\,\{\pm\epsilon_i\pm\epsilon_j,\pm2\epsilon_i\mid i,j\ge3\}\cr
\Phi=&\,\{\pm(\epsilon_1+\epsilon_1'+\epsilon_2')\}
\cup\{\pm(\epsilon_2+\epsilon_1'-\epsilon_2')\}
\cup\{\pm\epsilon_i\mid i\ge3\}
\cup\{\pm\epsilon_1',\pm\epsilon_2'\}\cr}
$$
Here, our algorithm terminates with
$\Phi_+^t=\{\epsilon_1+\epsilon_1'+\epsilon_2',
\epsilon_2+\epsilon_1'-\epsilon_2',\epsilon_1',\epsilon_2'\}$ and
$\Phi_0^s=\{\pm\epsilon_i\mid i\ge3\}$. Since $\Phi_+^t$ is linearly
independent, we see that $V$ is multiplicity free of rank $4$. The
Levi subgroup attached to $\Delta_0$ is $(\CC^*)^2\times
Sp_{2m-4}(\CC)\times(\CC^*)^2$. Since there is one equivalence class
of singular weights we see that the generic isotropy group is
$Sp_{2m-5}(\CC)$.

\medskip

The algorithm is very convenient in any given special
case, but it becomes quite awkward if large numbers of representations
(like series) have to be handled. Therefore, we prefer to use a couple
of ``shortcuts'' which we collect below. A very useful
criterion is:

\medskip

\ITEM{2}{A} {\it Let $(\fg,V)$ be a MFSR. Then
$\|dim|V\le\|dim|\fg+\|rk|\fg$.}

\medskip

\Proof: From $\|dim|\fg=|\Delta|+\|rk|\fg$ we get
$\|dim|V-\|dim|G+\|rk|G=|\Phi|-|\Delta|$. Now we run the
algorithm. From \cite{E5} we get
$|\Phi|-|\Delta|=|\Phi_0|-|\Delta_0|=2|\Phi_+^t|+|\Phi_0^s|-|\Delta_0|$.
The assertion follows from $|\Phi_+^t|\le\|rk|\fg$ (since $\Phi_+^t$
is linearly independent) and $|\Phi_0^s|\le|\Delta_0|$ (since
$2\Phi_0^s\subseteq\Delta_0$).\qed

\noindent The remaining criteria are dealing with the following
situation: $\fg=\fg_1+\fg_2+\fg_3$ and $V=V_1\oplus V_2$ such
that $\fg_1$ acts only on $V_1$ and $\fg_3$ acts only on $V_2$. Let
$\fl_1$ be the generic isotropy algebra of $\fg_1+\fg_2$ in $V_1$ and
denote its image in $\fg_2$ by $\fl_{12}$.

\medskip

\ITEM{5}{B} {\it Assume $\fl_1$ is reductive.Then $(\fg,V)$ is
multiplicity free if and only if $(\fg_1+\fg_2,V_1)$ and
$(\fl_{12}+\fg_3,V_2)$ are multiplicity free.}

\medskip

\Proof: This is a consequence of the algorithm if we apply it first
only to extremal weights of $V_1$ until all of them are toroidal and
then to the weights of $V_2$.\qed

\noindent The problem with criterion \cite{I5} is that
$(\fl_{12}+\fg_3,V_2)$ is in general not saturated. The most important
special case is:

\medskip

\ITEM{6}{C} {\it Assume $V_2=U_2\otimes U_3$ where $U_i$ is an
$\fg_i$\_module. If $\|dim|U_2\ge3$ and $\fl_{12}=0$ then $V$ is not
multiplicity free.}

\medskip

\Proof: This follows from \cite{I5} and the claim that in a MFSR no
highest weight $\chi$ appears more than twice. The claim is a
consequence of the algorithm since the multiplicity of $\chi$ goes
down at most once and then by at most one.\qed

\noindent
Another useful consequence of \cite{I5} is

\medskip

\ITEM{3}{D} {\it Assume
$$0
\eqalign{\hbox{either}\quad&
\fg_1=\sp_{2m},\ \fg_2=\so_n,\ V_1=\CC^{2m}\otimes\CC^n,\cr
\hbox{or}\quad&
\fg_1=\sl_m+\ft^1,\ \fg_2=\sl_n,\ V_1=T(\CC^m\otimes\CC^n).\cr}
$$
Assume that $V$ is not multiplicity free for $m=m_0\ge1$. Then $V$
is not multiplicity free for all $m\ge m_0$.}

\medskip

\Proof: Indeed $\fl_{12}(m)\subseteq\fl_{12}(m_0)$ for $m\ge m_0$ (see
Tables~1 and~2).\qed

\noindent The next criterion is useful in dealing with type~2
components:

\medskip

\ITEM{1}{E} {\it Assume $\fg_1=\sl_2$ and $V_1=\CC^2\otimes U$ where
$U$ is an $\fg_2$\_module. Suppose $(\fg,V)$ is not multiplicity
free. Then $(\ft^1+\fg_2+\fg_3,T(U)\oplus V_2)$ isn't either.}

\medskip

\Proof: Follows from $\ft^1\subseteq\sl_2$.\qed

\noindent Sometimes, we need a more refined argument:

\ITEM{8}{F} {\it Same setup as in \cite{I1} with $\|dim|U\ge2$. Let
$\fl$ be the generic isotropy algebra of $(\fg,V)$ and let $\fh$ be its
image in $\fg_1$. Then $(\ft^1+\fg_2+\fg_3,T(U)\oplus V_2)$ is
multiplicity free if and only if $(\fg,V)$ is multiplicity free and
$\fh\ne0$.}

\medskip

\Proof: Let $(\Phi_0,\Delta_0)$ and
$(\overline\Phi_0,\overline\Delta_0)$ be the outputs of the algorithm
for $(\fg,V)$ and $(\ft^1+\fg_2+\fg_3,T(U)\oplus V_2)$
respectively. Then $\overline\Delta_0=\Delta_0$ and
$\overline\Phi_0=\Phi_0\cup\{\chi-\alpha\}$ where $\chi$ is a highest
weight of $\CC^2\otimes U$ and $\alpha$ is the positive root of
$\sl_2$. Therefore, $\overline\Phi_+^t$ is linearly independent if and
only if $\Phi_+^t$ is linearly independent, i.e. $(\fg,V)$ is
multiplicity free, and $\alpha\not\in\<\Phi_+^t\>$, i.e. $\fh\ne0$.\qed

\beginsection Beweis. Proofs of the classification theorems

We start with a Lemma:

\Lemma. Let $\rho:\fg\pfeil\rho(V)$ be a MFSR. Then the centralizer of
$\rho(\fg)$ in $\sp(V)$ is commutative.

\Proof: Recall the description of the centralizer in the proof of
\cite{satcrit}. Thus, we have to show: the multiplicity of a component
of type~1, 2a, and 2b is $\le2$, $\le1$, and $\le0$,
respectively. This follows from the fact that the symplectic
representations $(\sp(V),V^{\oplus3})$, $(\gl(U),(U\oplus U^*)^{\oplus2})$ and
$(\so(U),U\oplus U)$ are not multiplicity free.\qed

\noindent{\it Proof of \cite{saturiert}:} The preceding lemma implies that
$\overline\fg$ is necessarily the normalizer of $\fg$ thereby proving
its uniqueness. Now assume $(\overline\fg,V)$ is multiplicity
free. Using the algorithm, we wind up with a linearly independent
subset $\Phi_+^t$ of $\overline\ft^*=\ft^*\oplus\fc^*$ where $\ft$ is
a Cartan subalgebra of $\overline\fg'$. The condition is now that
$\Phi_+^t$ stays linearly independent after restriction to
$\ft\oplus\fd$. The kernel of the restriction map is
$\fd^\perp\subseteq\fc^*$. Let
$S:=\<\Phi_+^t\>\subseteq\ft^*\oplus\fc^*$. Then the condition is
$S\cap\fd^\perp=0$. Since $S\cap\fd^\perp=(S\cap\fc^*)\cap\fd^\perp$
this is equivalent to $\fa+\fd=\fc$ with
$\fa=(S\cap\fc^*)^\perp\subseteq\fc$.\qed

Next we state the soundness of the tables:

\Lemma. All representations in Tables~1, 2, 11, 12, and 22 are
multiplicity free with the stated generic isotropy Lie algebra.

\Proof: We leave it to the reader to apply the algorithm on each
item. Note that \cite{TabS} below can be checked simultaneously
without extra effort.\qed

\noindent{\it Proof of \cite{IndecompClass}:} Here we rely on previous
classification work. Let $V$ be indecomposable of type~1. Thus, by
\cite{cofree}, $V$ is an irreducible cofree representation. These have
been classified by Littelmann in \cite{Lit}. If one extracts from the
tables in \cite{Lit} all representations which carry a symplectic
structure then one obtains exactly Table~1.

If $V=T(U)$ is of type~2 then $U$ is a multiplicity free space
(see \S\cite{Introduction}). The irreducible ones have been classified by
Kac~\cite{Kac} and are reproduced in Table~2.\qed

\noindent
During the rest of the classification we the tacitly use the following
observations:

\medskip

\item{$\bullet$} {\it Let $V_1$ and $V_2$ be two symplectic
representations. Then in order for $V=V_1\oplus V_2$ to be
multiplicity free it is necessary that both $V_1$ and $V_2$ are
multiplicity free.}  In fact, $\cW(V_i)^G$ is a subalgebra of
$\cW(V)^G$. 

\medskip

\item{$\bullet$} {\it Let $(\fg,V)$ be a symplectic representation and
$\fh\subseteq\fg$ a reductive subalgebra. Then in order for $(\fh,V)$
to be multiplicity free it is necessary that $(\fg,V)$ is
multiplicity free.} In fact, $\cW(V)^\fg$ is a subalgebra of
$\cW(V)^\fh$.

\medskip

\item{$\bullet$} {\it Each type~1 component appears with multiplicity
one}. In fact, otherwise the representation wouldn't be saturated
(\cite{satcrit}).

\medskip

\item{$\bullet$} {\it We may assume that $V$ contains at least one
type~1 component}. In fact, all representations consisting entirely of
type~2 components come from multiplicity free actions (see
\S\cite{Introduction}). Those have been classified by
Benson\_Ratcliff~\cite{BR} and Leahy~\cite{Leahy}. Their result is
reproduced in Table~22.

\medskip

\noindent
We continue by classifying all decomposable, connected, saturated
MFSRs where exactly one simple factor of $\fg$ acts on more than one
component effectively. This means
$$1
\fg=\fg_0+\fg_1+\ldots+\fg_s, V=V_1\oplus\ldots\oplus V_s,\ s\ge2
$$
where $\fg_0$ is simple and $V_i$ is an indecomposable symplectic
representation of $\fg_0+\fg_i$ for $i=1,\ldots,s$.
The main case will
be $s=2$ where we use the $\vplus$\_notation. For $s\ge3$ we use a
notation like, e.g.,
$$
\underline\so_n\otimes\sl_2\oplus\underline\so_n\otimes\sl_2
\oplus\underline\so_n\otimes\sl_2
$$
where $\fg_0$ is underlined. We actually show that $s=3$ is impossible
which settles also all higher cases.

\parindent0pt
\def\\{\itemitem{}}

\Lemma. There are no saturated MFSRs of the form \cite{E1} with
$\fg_0$ exceptional.

\Proof: The only components involving an exceptional algebra are
$$
\sl_2\otimes\sG2, \sp_4\otimes\sG2, T(\sG2), T(\sE6), \sE7
$$
We can exclude the cases $\sE7$, $T(\sE6)$,
and $\sp_4\otimes\sG2$ \cite{I6}. The leaves

\\$\sl_2\otimes\sG2\vplus\sG2\otimes\sl_2$: not multiplicity free
by~\cite{I2}.

\\$\sl_2\otimes\sG2\vplus T(\sG2)$: use~\cite{I1}.\qed

\medskip

\Lemma. All saturated MFSRs of the form \cite{E1} with
$\fg_0=\so_n$, $n\ge7$ are contained in Tables~11, 12, and 22.

\Proof: The only possible components not involving spin
representations are
$$2
\sp_{2m}\otimes\so_n\ (m\ge1),\quad T(\so_n)
$$

\\$\sp_{2m}\otimes\so_n\vplus\so_n\otimes\sp_{2p}$ ($m\ge p\ge1$): these
representations are multiplicity free for $m=p=1$. To exclude the
other cases we may assume $m=2$ and $p=1$~\cite{I3}. Then we use
\cite{I5}. We have $\fl_{12}=\so_{n-4}+\ft^2$. Thus
$(\fl_{12}+\sp_2,V_2)=\so_{n-4}\otimes\sl_2\vplus T(\sl_2)\vplus
T(\sl_2)$. Now we apply \cite{I6} (with $V_1$ equal to the last two
summands).

\\$\sp_{2m}\otimes\so_n\vplus T(\so_n)$ ($m\ge1$): not multiplicity
free by~\cite{I8}.
\medskip

The only way to combine the representations to a triple link (without
spin representation) is

\\$\underline\so_n\otimes\sl_2\oplus\underline\so_n\otimes\sl_2
\oplus\underline\so_n\otimes\sl_2$: we use \cite{I5} with $V_1$ being
the first two summands. Then $\fl_{12}=\so_{n-4}$. The representation
$(\fl_{12}+\fg_3,V_2)$ contains $(\sl_2,(\CC^2)^{\oplus4})$ which is not
multiplicity free by~\cite{I2}.

\medskip Now we consider spin representations. A glance at Tables~1
and~2 shows that only the cases $n\le13$ have to be
considered.

\medskip

$n=13$: In addition to \cite{E2} we have $\spin_{13}$.

\\$\sp_{2m}\otimes\so_{13}\vplus\spin_{13}$ ($m\ge1$): not multiplicity
free. First reduce to $m=1$~\cite{I3} and then use~\cite{I2}.

\\$T(\so_{13})\vplus\spin_{13}$: not multiplicity free~\cite{I1}.

\medskip

$n=12$: In addition to \cite{E2} we have both spin representations
$\spin_{12}^\pm$. Thus, up to isomorphism we have to check the
following cases:

\\$\sp_{2m}\otimes\so_{12}\vplus\spin_{12}$ ($m\ge1$): multiplicity free
for $m=1$ and $m=2$. For $m\ge3$ use~\cite{I3} and~\cite{I2}.

\\$T(\so_{12})\vplus\spin_{12}$: multiplicity free. 

\\$\spin_{12}^+\vplus\spin_{12}^-$: multiplicity free.

\medskip
For triple links we have the following possibilities:

\\$\underline\spin_{12}\oplus\underline\so_{12}\otimes\sl_2
\oplus\underline\so_{12}\otimes\sl_2$: we use \cite{I5} with $V_1$
being the last two summands. Then $\fl_{12}=\so_8$. The restriction of
$\spin_{12}$ to $\fl_{12}$ is $2\spin_8^++2\spin_8^-$, a
representation which is not multiplicity free (it would be though if
additionally $\ft^2$ were acting).

\\$\underline\spin_{12}^+\oplus\underline\spin_{12}^-
\oplus\underline\so_{12}\otimes\sp_{2m}$ ($m=1,2$): use~\cite{I2}.

\\$\underline\spin_{12}^+\oplus\underline\spin_{12}^-
\oplus T(\underline\so_{12})$: use~\cite{I1}.

\medskip

$n=11$: In addition to \cite{E2} we have $\spin_{11}$. 

\\$\sp_{2m}\otimes\so_{11}\vplus\spin_{11}$ ($m\ge1$): multiplicity free
for $m=1$. For $m\ge2$ use~\cite{I3} and~\cite{I2}.

\\$T(\so_{11})\vplus\spin_{11}$: not multiplicity free~\cite{I8}.

Possible triple links:

\\$\underline\spin_{11}\oplus\underline\so_{11}\otimes\sl_2
\oplus\underline\so_{11}\otimes\sl_2$: use~\cite{I2}.

\medskip

$n=10$: In addition to \cite{E2} we have $T(\spin_{10})$ (observe
$T(\spin_{10}^+)=T(\spin_{10}^-)$). 

\\$\sp_{2m}\otimes\so_{10}\vplus T(\spin_{10})$ ($m\ge1$): multiplicity free
for $m=1$. For $m\ge2$ use~\cite{I3} and~\cite{I2}.

Possible triple links are:

\\$T(\underline\spin_{10})\oplus\underline\so_{10}\otimes\sl_2
\oplus\underline\so_{10}\otimes\sl_2$: use~\cite{I2}.

\medskip

$n=9$: In addition to \cite{E2} we have $\sl_2\otimes\spin_9$ and
$T(\spin_9)$. 

\\$\sp_{2m}\otimes\so_9\vplus\spin_9\otimes\sl_2$ ($m\ge1$): not
multiplicity free~\cite{I3},~\cite{I2}.

\\$\sp_{2m}\otimes\so_9\vplus T(\spin_9)$ ($m\ge1$): not multiplicity
free~\cite{I1}.

\\$\sl_2\otimes\spin_9\vplus\spin_9\otimes\sl_2$: not multiplicity free~\cite{I2}.

\\$\sl_2\otimes\spin_9\vplus T(\spin_9)$: not multiplicity
free~\cite{I1}.

\\$\sl_2\otimes\spin_9\vplus T(\so_9)$: not multiplicity
free~\cite{I1}.

\medskip

$n=8$: Formally, there are no new representations than \cite{E2} but
due to triality we may replace $\so_8$ by either spin
representation. Thus, up to isomorphism we encounter the following
cases:

\\$\sp_{2m}\otimes\so_8\vplus\spin_8\otimes\sp_{2p}$ ($m\ge p\ge1$):
multiplicity free for $(m,p)=(1,1)$ and $(2,1)$. Using~\cite{I3} it
remains to check $(m,p)=(2,2)$ and $(3,1)$. Both cases are handled
by~\cite{I2}.

\\$\sp_{2m}\otimes\so_8\vplus T(\spin_8)$: multiplicity free only for
$m=1$. Use~\cite{I8} for $m\ge2$.

\medskip

For triple links we may ignore $T(\so_8)$ by~\cite{I1}. Thus, we are
left with

\\$\rho_1(\underline\so_8)\otimes\sp_{2m}
\oplus\rho_2(\underline\so_8)\otimes\sl_2
\oplus\rho_3(\underline\so_8)\otimes\sl_2$ where $m=1,2$ and $\rho_i$
denotes either of the three 8\_dimensional fundamental representations. In all
these cases we can use~\cite{I2}.

\medskip

$n=7$: In addition to \cite{E2} we have $\sp_{2m}\otimes\spin_7$
($m\ge1$) and $\ T(\spin_7)$. For $m\ge3$, the image of the generic
isotropy algebra of $\sp_{2m}\otimes\spin_7$ in $\spin_7$ is
trivial. Therefore, we may assume $m\le2$ \cite{I6}.

\\$\sp_{2m}\otimes\so_7\vplus\spin_7\otimes\sp_{2p}$ ($m,p\ge1$):
multiplicity free for $m=p=1$. Using~\cite{I3} and the remark above
we have to check $(m,p)=(2,1)$ and $(1,2)$. To do this use~\cite{I2}.

\\$\sp_{2m}\otimes\spin_7\vplus\spin_7\otimes\sp_{2p}$ ($m,p\ge1$):
multiplicity free for $m=p=1$. It remains to check the three cases
$(m,p)=(2,1)$, $(1,2)$, and $(2,2)$ using~\cite{I2}.

\\$\sp_{2m}\otimes\so_7\vplus T(\spin_7)$ ($m\ge1$): multiplicity
free only for $m=1$. For $m\ge2$ use~\cite{I1}.

\\$\sp_{2m}\otimes\spin_7\vplus T(\spin_7)$: not multiplicity
free~\cite{I8}.

\\$\sp_{2m}\otimes\spin_7\vplus T(\so_7)$: not multiplicity
free~\cite{I8}.

\medskip
For triple links we may ignore $T(\so_7)$ and
$T(\spin_7)$ by~\cite{I1}. Then we are left with the representations

\\$\rho_1(\underline\so_7)\otimes\sl_2
\oplus\rho_2(\underline\so_7)\otimes\sl_2
\oplus\rho_3(\underline\so_7)\otimes\sl_2$ where $\rho_i$ denotes
either the defining or the spin representation. All cases are ruled
out with~\cite{I2}.\qed

\Lemma. All saturated MFSRs of the form \cite{E1} with
$\fg_0=\sl_n$, $n\ge3$ are contained in Tables~11, 12, and 22.

\Proof: Indecomposable MFSRs of type~1 involve only $\sl_6$ and
$\sl_4=\so_6$.

\medskip

$n=6$: We have to consider
$$
\wedge^3\sl_6, T(\sl_6\otimes\sl_m)\ (m\ge2),
T(\wedge^2\sl_6), T(S^2\sl_6), T(\sl_6), T(\sp_4\otimes\sl_6).
$$

\\$\wedge^3\sl_6\vplus T(\sl_6\otimes\sl_m)$ ($m\ge2$): multiplicity free
for $m=2$. For $m\ge3$ use~\cite{I3} and~\cite{I2}.

\\$\wedge^3\sl_6\vplus T(\wedge^2\sl_6)$: use~\cite{I2}.

\\$\wedge^3\sl_6\vplus T(S^2\sl_6)$: use~\cite{I2}.

\\$\wedge^3\sl_6\vplus T(\sl_6)$: multiplicity free.

\\$\wedge^3\sl_6\vplus T(\sl_6\otimes\sp_4)$: use~\cite{I2}.

\medskip
The only possible representations with triple link are:

\\$\wedge^3\underline\sl_6\oplus T(\underline\sl_6)
\oplus T(\underline\sl_6)$: This is a close call. The algorithm gives
7 weights which turn out to be linearly dependent (see the first
example illustrating the algorithm): not multiplicity free.

\\$\wedge^3\underline\sl_6\oplus T(\underline\sl_6)
\oplus T(\underline\sl_6\otimes\sl_2)$: use~\cite{I2}.

\medskip

$n=4$: Observe $\sl_4=\spin_6$ and $\wedge^2\sl_4=\so_6$. Then we get the
representations
$$
\sp_{2m}\otimes\so_6,
T(\sl_m\otimes\sl_4)\ (m\ge2),
T(\so_6),
T(S^2\sl_4),
T(\sl_4),
T(\sp_4\otimes\sl_4).
$$

\\$\sp_{2m}\otimes\so_6\vplus\so_6\otimes\sp_{2p}$ ($m,p\ge1$): same
argument as for general $\so_n$: multiplicity free only for $m=p=1$.

\\$\sp_{2m}\otimes\so_6\vplus T(\so_6)$ ($m\ge1$): use~\cite{I8}.

\\$\sp_{2m}\otimes\so_6\vplus T(\sl_4\otimes\sl_p)$ ($m\ge1$,
$p\ge2$): multiplicity free for $(m,p)=(1,2)$. Using~\cite{I3}
we can reduce the remaining cases to $(m,p)=(2,2)$ and
$(1,3)$. Use~\cite{I2} for these cases.

\\$\sp_{2m}\otimes\so_6\vplus T(S^2\sl_4)$: the generic isotropy group
of $T(S^2\sl_4)$ is trivial. Now use~\cite{I6}.

\\$\sp_{2m}\otimes\so_6\vplus T(\sl_4)$: multiplicity free for all $m\ge1$

\\$\sp_{2m}\otimes\so_6\vplus T(\sl_4\otimes\sp_4)$: use~\cite{I3}
and~\cite{I2}.

\medskip

The possible representation with triple link are:

\\$\underline\so_6\otimes\sl_2\oplus
\underline\so_6\otimes\sl_2\oplus\underline\so_6\otimes\sl_2$: use~\cite{I2}.

\\$\underline\so_6\otimes\sl_2\oplus
\underline\so_6\otimes\sl_2\oplus T(\underline\sl_4)$: use~\cite{I2}.

\\$\underline\so_6\otimes\sl_2\oplus
\underline\so_6\otimes\sl_2\oplus T(\underline\sl_4\otimes\sl_2)$: use~\cite{I2}.

\\$\underline\so_6\otimes\sp_{2m}\oplus T(\underline\sl_4)
\oplus T(\underline\sl_4)$ ($m\ge1$): first reduce to $m=1$
with~\cite{I3}. Then use~\cite{I2}.

\\$\underline\so_6\otimes\sl_2\oplus T(\underline\sl_4)
\oplus T(\underline\sl_4\otimes\sl_2)$: use~\cite{I2}.\qed

\Lemma. All saturated MFSRs of the form \cite{E1} with
$\fg_0=\sp_{2n}$, $n\ge2$ are contained in Tables~11, 12, and 22.

\Proof: For general $n$, we have to consider the components
$$
\so_m\otimes\sp_{2n}\ (m\ge3),\ \sp_{2n}\otimes\spin_7,\
\sp_{2n},\ T(\sp_{2n}),\ T(\sp_{2n}\otimes\sl_2),\
T(\sp_{2n}\otimes\sl_3).
$$
This leaves the following possibilities with $s=2$:

\\$\so_m\otimes\sp_{2n}\vplus\sp_{2n}$ ($m\ge3$): multiplicity free for
all $m\ge3$, $n\ge2$.

\\$\so_m\otimes\sp_{2n}\vplus\sp_{2n}\otimes\so_p$ ($m\ge p\ge3$): here,
we directly apply the algorithm. For $m=3$ and $p\ge2$ we get the
linearly dependent chain of weights
$$0
\epsilon_1+\epsilon'_1,\epsilon'_1+\epsilon''_1,
\epsilon'_1-\epsilon''_1,\epsilon'_2+\epsilon''_1,\epsilon'_2.
$$
If $m\ge4$ and $p\ge2$ we get
$$0
\epsilon_1+\epsilon'_1,\epsilon_2+\epsilon'_2,
\epsilon'_1+\epsilon''_1,\epsilon'_1-\epsilon''_1,
\epsilon'_2+\epsilon''_1,\epsilon'_2-\epsilon''_1
$$
which is also linearly dependent.

Observe that this argument works
also for $p=2$. Thus we get:

\\$\so_m\otimes\sp_{2n}\vplus T(\sp_{2n})$ ($m\ge3$): not multiplicity free

\\$\so_m\otimes\sp_{2n}\vplus\sp_{2n}\otimes\spin_7$ ($m\ge3$): use
$\spin_7\subset\so_8$.

\\$\so_m\otimes\sp_{2n}\vplus T(\sp_{2n}\otimes\sl_2)$ ($m\ge3$): we
use \cite{I5} with $V_1$ being the second summand. Then
$\fl_{12}=\sp_{2n-4}+\ft^1$. This implies that $(\fl_{12}+\so_m,V_2)$
contains four summands of $\so_m$. Its saturation, $T(\so_m)\vplus
T(\so_m)$, is not multiplicity free.

\\$\so_m\otimes\sp_{2n}\vplus T(\sp_{2n}\otimes\sl_3)$ ($m\ge3$): same
argument as above.

\\$\sp_{2n}\vplus T(\sp_{2n}\otimes\sl_2)$: use a refinement of the
argument. This time a single $\ft^1$ is acting on $\CC^4$.

\\$\sp_{2n}\vplus T(\sp_{2n}\otimes\sl_3)$: same
argument as above.

\\$\sp_{2n}\vplus\sp_{2n}\otimes\spin_7$: multiplicity free for
$n=2$. For $m\ge3$ use~\cite{I6} on the second summand.

\\$\sp_{2n}\vplus T(\sp_{2n})$: multiplicity free for all $n\ge1$.

\\$\spin_7\otimes\sp_{2n}\vplus\sp_{2n}\otimes\spin_7$: use
$\spin_7\subset\so_8$. 

\\$\spin_7\otimes\sp_{2n}\vplus T(\sp_{2n})$: use
$\spin_7\subset\so_8$. 

\\$\spin_7\otimes\sp_{2n}\vplus T(\sp_{2n}\otimes\sl_2)$: use
$\spin_7\subset\so_8$. 

\\$\spin_7\otimes\sp_{2n}\vplus T(\sp_{2n}\otimes\sl_3)$: use
$\spin_7\subset\so_8$. 

\medskip

The possible representation with $s=3$ are so far:

\\$\underline\sp_{2n}\oplus T(\underline\sp_{2n})\oplus
T(\underline\sp_{2n})$: we use \cite{I5} with $V_1$ being the last two
summands. Thus $\fl_{12}=\sp_{2n-4}$. It has a 4\_dimensional fixed
point space on the first summand: not multiplicity free.

\medskip

$n=6$: we have to consider additionally $\wedge^3_0\sp_6$.

\\$\so_m\otimes\sp_6\vplus\wedge^3_0\sp_6$: we use \cite{I5} with
$V_1$ the second summand. Then $\fl_{12}=\sl_3$. The saturation of
$(\so_m+\fl_{12},V_2)$ is $T(\so_m\otimes\sl_3)$ which is not
multiplicity free (since it is not in Table~2).

\\$\spin_7\otimes\sp_6\vplus\wedge^3_0\sp_6$: use
$\spin_7\subset\so_8$.

\\$\sp_6\vplus\wedge^3_0\sp_6$: multiplicity free.

\\$T(\sp_6)\vplus\wedge^3_0\sp_6$: multiplicity free.

\\$T(\sl_2\otimes\sp_6)\vplus\wedge^3_0\sp_6$: use~\cite{I2}.

\\$T(\sl_3\otimes\sp_6)\vplus\wedge^3_0\sp_6$: use~\cite{I2}.

\medskip

The possible representations with triple link are:

\\$\wedge^3_0\underline\sp_6\oplus\underline\sp_6
\oplus T(\underline\sp_6)$: use~\cite{I2}.

\\$\wedge^3_0\underline\sp_6\oplus T(\underline\sp_6)
\oplus T(\underline\sp_6)$: use~\cite{I2}.

\medskip

$n=4$: Observe $\sp_4=\spin_5$. Then we get additionally
$\sp_{2m}\otimes\so_5$, $\sp_4\otimes\sG2$, $T(\sp_4\otimes\sl_m)$
($m\ge2$), $T(\so_5)$

\\$\sp_{2m}\otimes\so_5\vplus\sp_4\otimes\so_p$
($m\ge1,p\ge3$): Using~\cite{I3} we may assume $m=1$. We use \cite{I5}
with $V_1$ being the first summand. Then $\fl_{12}=\so_3+\ft^1$. The
restriction of $V_2$ to $\fl_{12}+\so_p$ is
$T(\sl_2\otimes\so_p)$ which is not multiplicity free.

\\$\sp_{2m}\otimes\so_5\vplus\sp_4\otimes\spin_7$ ($m\ge1$): use
$\spin_7\subset\so_8$.

\\$\sp_{2m}\otimes\so_5\vplus\sp_4$ ($m\ge1$): multiplicity free.

\\$\sp_{2m}\otimes\so_5\vplus T(\sp_4)$ ($m\ge1$): multiplicity free
for $m=1$. For the other cases use~\cite{I3} and then~\cite{I2}.

\\$\sp_{2m}\otimes\so_5\vplus \so_5\otimes\sp_{2p}$ ($m\ge p\ge1$):
multiplicity free for $m=p=1$. The argument for general $\so_n$ works
here as well.

\\$\sp_{2m}\otimes\so_5\vplus T(\sp_4\otimes\sl_p)$ ($m\ge1,p\ge2$):
first, we reduce to $m=1$ \cite{I3}. Then, as in the first case, we obtain
$\sl_2+\sl_p+\ft^2$ acting on $T(\CC^2\otimes\CC^p)\oplus
T(\CC^2\otimes\CC^p)$ which is not multiplicity free.

\\$\sp_{2m}\otimes\so_5\vplus T(\so_5)$ ($m\ge1$): not multiplicity
free by~\cite{I8}.

\\$\so_m\otimes\sp_4\vplus T(\sp_4\otimes\sl_p)$ ($m\ge3,p\ge2$): for
$p\ge3$, the image of $\fl(V_2)$ in $\sp_4$ is $0$. The case $p=2$ has
already been dealt with.

\\$\spin_7\otimes\sp_4\vplus T(\sp_4\otimes\sl_p)$ ($p\ge2$): use
$\spin_7\subset\so_8$.

\\$\sp_4\vplus T(\sp_4\otimes\sl_p)$: for $p\ge3$ argue as above. For
$p=2$ use~\cite{I2}.

\\$\so_m\otimes\sp_4\vplus T(\so_5)$ ($m\ge3$): use~\cite{I1}.

\\$\spin_7\otimes\sp_4\vplus T(\so_5)$: use $\spin_7\subset\so_8$

\\$\sp_4\vplus T(\so_5)$: multiplicity free.

\medskip
The possible representation with triple link are:

\\$\underline\so_5\otimes\sl_2\oplus\underline\so_5\otimes\sl_2
\oplus\underline\so_5\otimes\sl_2$,

\\$\underline\so_5\otimes\sl_2\oplus\underline\so_5\otimes\sl_2
\oplus\underline\sp_4$,

\\$\underline\so_5\otimes\sl_2\oplus\underline\so_5\otimes\sl_2
\oplus T(\underline\sp_4)$,

\\$\underline\so_5\otimes\sl_2
\oplus \underline\sp_4\oplus T(\underline\sp_4)$,

\\$\underline\so_5\otimes\sl_2
\oplus T(\underline\sp_4)\oplus T(\underline\sp_4)$, and

\\$\underline\sp_4\oplus T(\underline\sp_4)\oplus T(\underline\sp_4)$:
all of them can be handled with~\cite{I2}.\qed

\parindent20pt

\noindent{\it Proof of \cite{CompClass}:} First, assume that $V$ has
only two components. If more than one simple factor $\not\cong\sl_2$
acts effectively on both components then it must be obtained by
identifying simple factors of $\fg_1$ with simple
factors of $\fg_2$. Going through Tables~11, 12, and 22 one sees that
this is impossible without creating an $\sl_2$\_link. Thus, these
tables contain in fact all two\_component, connected, saturated MFSRs.

If $V$ has three components then we have already seen that no simple
factor of $\fg$ can act effectively on all of them. The other scenario
is
$$
\fg=\fg_0'+\fg_0''+\fg_1+\fg_2+\fg_3, V=V_1\oplus V_2\oplus V_3
$$
where $\fg_0'$ and $\fg_0''$ are simple, $V_1$ is an
$\fg_1+\fg_0'$\_module, $V_2$ is an $\fg_0'+\fg_2+\fg_0''$\_module and
$V_3$ is an $\fg_0''+\fg_3$\_module. I claim that no multiplicity free
representations of that form exist.

Both representations $(\fg_0'+\fg_0''+\fg_1+\fg_2, V_1\oplus
V_2)$ and $(\fg_0'+\fg_0''+\fg_2+\fg_3,V_2\oplus V_3)$ have to occur in
one of the Tables~11, 12, or 22. Here are the only possibilities to combine
two of them:

\\$\spin_{12}\vplus\so_{12}\otimes\sp_4\vplus\sp_4$,

\\$\sl_2\otimes\spin_8\vplus\so_8\otimes\sp_4\vplus\sp_4$,

\\$T(\sl_4)\vplus\so_6\otimes\sp_{2m}\vplus\sp_{2m}$ $(m\ge2$), and

\\$\sp_4\vplus\so_5\otimes\sp_{2m}\vplus\sp_{2m}$ $(m\ge2$): in all
cases use~\cite{I6} on the first two summands.

\noindent Thus, all possibilities are exhausted and \cite{CompClass} is
proven.\qed

Before we start with the proof of \cite{sl2link}, we state the
significance of Table~S.

\Lemma TabS. Let $(\fg,V)$ be a saturated MFSR without $\sl_2$\_link
and $\Omega$ a non\_empty subset of the set of $\sl_2$\_factors of
$\fg$ acting non\_trivially on $V$. Let $\Phi_+^t$ be ``the'' output
of the algorithm applied to $(\fg,V)$ and let $A$ be the set of simple
roots of the $sl_2$\_factors in $\Omega$. Assume $\Phi_+^t\cup A$ is
linearly independent and that $\Omega$ is maximal with that
property. Then $(V,\Omega)$ is equivalent to an entry of Table~S where
we indicate membership in $\Omega$ by underlining.

\Proof: This is easily checked using the algorithm. Observe that the
condition to be checked is independent of the choices during the
algorithm since it can be rephrased as: let $\tilde\fg\subset\fg$ be
the subalgebra obtained by replacing each $\sl_2$\_factor in $\Omega$
by its Cartan subalgebra. Then $(\tilde\fg,V)$ is still multiplicity
free (see the proof of criterion \cite{I8}).\qed

\noindent{\it Proof of \cite{sl2link}:} Let $(\fg,V)$ be a saturated
symplectic representation and let $\fs\subseteq\fg$ be an
$\sl_2$\_link. Let $V=V_1\oplus\ldots\oplus V_r$ be the decomposition
of $V$ into components. Assume $\fs$ acts non\_trivially precisely on
$V_1,\ldots,V_s$ with $2\le s\le r$. Let $\overline\fg$ equal $\fg$
but with $\fs$ replaced by $\fs_1+\ldots+\fs_s$ where each
$\fs_i\cong\sl_2$. Then $(\overline\fg,V)$ is still a saturated
symplectic representation and none of the $\fs_i$ is an
$\sl_2$\_link. Conversely, $\fg$ is obtained from $\overline\fg$ by
replacing $\fs_1+\ldots+\fs_s\subseteq\overline\fg$ by its
diagonal. Observe that, by construction, each $\fs_i$ acts on a
different component of $V$ (this explains the first exception of the
theorem). We want to find a criterion for when $(\fg,V)$ is
multiplicity free. Clearly, $(\overline\fg,V)$ has to be multiplicity
free to begin with.

For $i=1,\ldots,s$ let $\alpha_i$ be the simple positive root of
$\fs_i$ and $\chi_i$ be a highest weight in $V_i$. First of all, at
most one component $V_i$ can be of type (S.9) since otherwise
$(\fg,V)$ wouldn't be saturated (this explains the second exception of
the theorem). Thus, if $V_1$ is not of type (S.9) then $\chi_1$ is an
extremal weight which is neither toroidal nor singular. Now we apply
our algorithm first to $\chi_1$ and then $\chi_2,\ldots,\chi_s$. If we
do that with respect to the root system of $\overline\fg$ then we have
to delete, among others, $\chi_1-\alpha_1,\ldots,\chi_s-\alpha_s$. If
we do it with respect to the root system of $\fg$ then we keep
$\chi_2-\alpha_2,\ldots,\chi_s-\alpha_s$ and we have to identify
$\alpha_1=\ldots=\alpha_s$. Since $\chi_1,\ldots,\chi_s$ are being
kept in any case we see that $s$ can't be bigger than 2. This shows that
that the identifications have to be pairwise and disjoint.

Thus, we have $s=2$. Then $\Phi_+^t$ (the output with respect to
$\fg$) differs from $\overline\Phi_+^t$ (the output with respect to
$\overline\fg$) by the additional weight $\chi_2-\alpha_2$ and the
identification $\alpha_1=\alpha_2$. This implies that $(\fg,V)$ is
multiplicity free if and only if
$\overline\Phi_+^t\cup\{\alpha_1,\alpha_2\}$ is linearly independent.

Since $s=2$, a general saturated MFSR is obtained from a collection of
representation from tables~1 through~22 by pairwise identifying
various $\sl_2$'s. Let $(\tilde\fg,\tilde V)$ one of these items and
let $\fs_1,\ldots,\fs_s$ be those $\sl_2$\_factors of $\tilde\fg$
which are being identified with other $\sl_2$\_factors. Then an
identification pattern is permissible if and only the union of
$\Phi_+^t$ with the roots of the $\fs_i$ is linearly independent.
This is precisely the condition for being member of Table~S
(\cite{TabS}).\qed

\beginrefs

\L|Abk:BR|Sig:BR|Au:Benson, C.; Ratcliff, G.|Tit:A classification of
multiplicity free actions|Zs:J. Algebra|Bd:181|S:152--186|J:1996|xxx:-||

\L|Abk:Howe|Sig:Ho|Au:Howe, R.|Tit:Remarks on classical invariant
theory|Zs:Trans. Amer. Math. Soc.|Bd:313|S:539--570|J:1989|xxx:-||

\L|Abk:HU|Sig:HU|Au:Howe, R.; Umeda, T.|Tit:The Capelli identity, the
 double commutant theorem, and multiplicity-free
 actions|Zs:Math. Ann.|Bd:290|S:565--619|J:1991|xxx:-||

\L|Abk:Kac|Sig:Ka|Au:Kac, V.|Tit:Some remarks on nilpotent orbits|Zs:J.
Algebra|Bd:64|S:190--213|J:1980|xxx:-||

\L|Abk:inprep|Sig:Kn|Au:Knop, F.|Tit:Invariant functions on
symplectic representations|Zs:Preprint|Bd:-|S:24
pages|J:2005|xxx:math.AG/0506171||

\L|Abk:Leahy|Sig:Le|Au:Leahy, A.|Tit:A classification of multiplicity free
representations|Zs:J. Lie Theory|Bd:8|S:367--391|J:1998|xxx:-||

\L|Abk:Lit|Sig:Li|Au:Littelmann, P.|Tit:Koregul\"are und
\"aquidimensionale Darstellungen|Zs:J. Algebra|Bd:123|S:193--222|J:1989|xxx:-||

\endrefs

\bye